\documentclass[11pt]{article}
\usepackage{latexsym,amssymb,amsmath,amsfonts,amsthm}
\usepackage{graphics,graphicx,mathrsfs,subfigure}
\usepackage{color,xspace}
\newcommand{\R}{{\mat R}}

\newcommand{\N}{{\mat N}}
\newcommand{\C}{{\mat C}}

\newcommand{\be}{\begin{eqnarray}}
\newcommand{\ben}{\begin{eqnarray*}}
\newcommand{\en}{\end{eqnarray}}
\newcommand{\enn}{\end{eqnarray*}}

\newcommand{\mat}{\mathbb}

\newtheorem{theorem}{Theorem}[section]
\newtheorem{lemma}[theorem]{Lemma}

\newtheorem{remark}[theorem]{Remark}
\newtheorem{proposition}[theorem]{Proposition}

\definecolor{rot}{rgb}{1,0,0}
\definecolor{hw}{rgb}{0,0,1}

\begin{document}
\renewcommand{\theequation}{\arabic{section}.\arabic{equation}}
%\begin{titlepage}
\title{\bf
Increasing stability for inverse acoustic source problems in the time domain
}
\author{Chun Liu\thanks{School of Mathematical Sciences, Nankai University,
Tianjing 300071, China  ({\tt liuchun@nankai.edu.cn})}
\and Suliang Si\thanks{School of Mathematics and Statistics, Shandong University of Technology,
Shandong 255000, China ({\tt sisuliang@amss.ac.cn})}
\and Guanghui Hu\thanks{School of Mathematical Sciences, Nankai University,
Tianjing 300071, China ({\tt ghhu@nankai.edu.cn})}
\and Bo Zhang\thanks{LSEC and Academy of Mathematics and Systems Sciences, Chinese Academy of Sciences,
Beijing, 100190, China and School of Mathematical Sciences, University of Chinese Academy of Sciences,
Beijing 100049, China ({\tt b.zhang@amt.ac.cn})}}
\date{}
%\end{titlepage}

%\vspace{.2in}

%\begin{document}
%\renewcommand{\theequation}{\arabic{section}.\arabic{equation}}

\maketitle

\begin{abstract}
This paper is concerned with inverse source problems for the acoustic wave equation in the full space $\R^3$, where the source term is compactly supported in both time and spatial variables. The main goal is to investigate increasing stability for the wave equation in terms of the interval length of given parameters (e.g., bandwith of the temporal component of the source function). We establish increasing stability estimates of the $L^2$-norm of the source function by using only the Dirichlet boundary data. Our method relies on the Huygens’ principle, the Fourier transform and explicit bounds for the continuation of analytic functions.
\end{abstract}

\section{Introduction}
Consider an acoustic wave propagation problem caused by a compactly supported source function in three dimensions. This can be modeled by the time-dependent wave equation
\begin{equation}\label{eq1}
  \partial_t^2u(x,t) -\lambda\Delta u(x,t)=F(x,t),\qquad  x\in\mathbb{R}^3,\ t>0,
\end{equation}
where $u(x,t)$ denotes the wave field, $0<\lambda\in\R$ and $F(x,t)$ is the source term. Together with the above governing equation, we impose the homogeneous initial conditions
\begin{equation}\label{eq2}
   u(x,0)=\partial_t u(x,0)=0,\qquad x\in\mathbb{R}^3.
\end{equation}
If $F(x,t)\in L^2([0,\infty);L^2(\R^3))$ has compact support, the problem (\ref{eq1})-(\ref{eq2}) admits a unique solution (see e.g., \cite{HW,LM,LM1})
\begin{equation}
 \nonumber
  u(x,t)\in{\mathcal{C}([0,+\infty);H^1(\mathbb{R}^3))\cap \mathcal{C}^1([0,+\infty);L^2(\mathbb{R}^3))}.
 \end{equation}
In this paper we consider three inverse problems depending on the a priori  knowledge on the parameter $\lambda$ and the form of the source function $F$. 

\textit{\textbf{IP1}}:
Assume $\lambda=1$ and $F(x,t)=f(x)g(t)$, where the temporal function $g$ is given. Suppose that $f,\ g$ have compact supports such that supp $f\subset B_R$ and supp $g\subset(0,T_0)$, where $ B_R:=\{x\in\mathbb{R}^3| \ |x|<R\}$ for some $R$ and $T_0>0$. The inverse problem is to recover the source term $f(x)$ from the Dirichlet boundary data $u(x,t)$ measured on $\partial B_R\times(0,T)$ with $T>0$ sufficiently large.

If we fix the parameter $\lambda$, it is in general impossible to uniquely recover a source term of the form $F(x,t)$, due to the presence of time-dependent non-radiating sources (see \cite{Bleistein1977Nonuniqueness,Bao2010multi-frequency}). 
Motivated by \cite{HKZ}, we suppose that the measurement data are given by 
a family of parameter-dependent functions $u_\lambda(x,t)$ for all 
$\lambda\in(0,\Lambda^2)$ with  $1<\Lambda<\infty$. Here $u_\lambda(x,t)$ denotes the unique solution to (\ref{eq1})-(\ref{eq2}) corresponding to the parameter $\lambda$.

\textit{\textbf{IP2}}:
Assume the source function $F(x,t)$ has compact supports such that supp$F(x,t)\subset B_R\times(0,T_0)$. The inverse problem is to determine  source term $F(x,t)$ from the boundary observation data $\{u_\lambda(x,t)|\  x\in\partial B_R, \ t\in(0,T), \ \lambda\in(0,\Lambda^2)\}$.

In the last inverse problem, we suppose that the $x_3$-dependent component of the source function is given.

\textit{\textbf{IP3}}:
Assume  $\lambda=1$ and $F(x,t)=f(\tilde{x},t)g(x_3)$, where the function $f$ satisfies supp$f\subset\widetilde{B}_{R_0}\times(0,T_0)$ and $g$ is supported in $(-R_0,R_0)$. Here $\tilde{x}=(x_1,x_2)\in\R^2$ and $\widetilde{B}_{R_0}:=\{\tilde{x}\in\R^2|\ |\tilde{x}|<R_0\}$  The inverse problem is to determine the source term $f(\tilde{x},t)$ from the boundary observation data $u(x,t)$ on $\partial B_R\times(0,T)$ with $R>0$, $T>0$ sufficiently large.

Inverse source problems have many significant applications in scientific and engineering areas, as an important research subject
in inverse scattering theory.
For instance, detection of submarines and non-destructive measurement of industrial objects can be regarded as recovery of acoustic sources from boundary measurements of the pressure. Other applications include biomedical imaging optical tomography \cite{Arridge1999Optical,Isakov1990inverse} and geophysics.  
Consequently, inverse source problems have continuously attracted much attention by many researchers \cite{Bao2010multi-frequency,BaoStability2020,BaoLu,Isakov1990inverse,Yamamoto1999Stability,Yamamoto1999Uniqueness} and a great deal of mathematical and numerical results are available, especially for the time-dependent and time-harmonic acoustic waves. 

Inverse source problems in the time domain are usually treated as hyperbolic systems by applying the Carleman estimates \cite{Klibanov1999Inverse} and unique continuation theory. The approach of Carleman estimates can be used to recover both 
coefficients and source functions for hyperbolic equations; we refer to \cite{Choulli2006Some, Jiang2016Inverse,Yamamoto1999Stability,Yamamoto1999Uniqueness} for an incomplete list . In the time-harmonic regime, it is well known that there is no uniqueness for inverse source problems with a single frequency due to the existence of non-radiating sources \cite{Bleistein1977Nonuniqueness,Hauer2005On}. Computationally, a more serious
issue is the lack of stability, i.e., a small variation of the data might lead to a huge error in the
reconstruction. Hence it is crucial to study the stability of inverse source problems. The use of multiple-frequency data is an effective way to overcome non-uniqueness and has received a lot of attention in recent years. The paper \cite{Eller2009Acoustic} show uniqueness and numerical results for the Helmholtz equation with multi-frequency data. In \cite{Bao2010multi-frequency}, Bao et al. firstly get increasing stability for the Helmholtz equation by applying direct spatial Fourier analysis methods. A different method is proposed 
in \cite{Cheng2016Increasing} to derive increasing stability bounds for the  Helmholtz equation in three dimensions, which involves
 the temporal
Fourier transform together with sharp bounds of the analytic continuation at higher wave numbers. The authors of \cite{Cheng2016Increasing} firstly bridge the Helmholtz equation in the frequency-domain with the associated hyperbolic equations to get  increasing stability results. These results were later generalized to the Helmholtz equation and Maxwell's equation in three dimensions (see \cite{Li2017Increasing} and \cite{BaoStability2020}). We also refer to \cite{HKZ,Z2023} for uniqueness
results of time-dependent inverse source problems using the approach of Fourier-Laplace transform. 

To the best of our knowledge, increasing stability results are obtained for the Helmholtz equation in the time-harmonic domain only. Even the concept of increasing stability seems not available in the time domain.
In this paper we are interested in the increasing stability of inverse source problems for the time-dependent acoustic equation in three dimensions. Inspired by the approach of Fourier transform, we first reduce inverse source problems of the acoustic equation to that of the associated Helmholtz equation with multiple frequencies. Then we derive increasing stability estimate by using sharp bounds of analytic continuation given in \cite{Cheng2016Increasing}. Instead of the Cauchy data used in the works mentioned above, we use only Dirichlet data on the lateral boundary. This work
initializes mathematical studies on increasing stability of inverse source problems
for the wave equation. %Our objective is to develop an increasing stability theory of the inverse
%source problems for the wave equation. 
The existing approaches dealing with the Helmholtz equations has been adapted to handle inverse source problems for the wave equation in the time domain.

The rest of this paper is organized as follows. In Sections \ref{man} and \ref{pre}, we state main results and well-posedness of the direct problem. Section \ref{th1} is devoted to the increasing stability of the inverse problem IP1. In Section \ref{th2}, we investigate the second problem IP2 for general source terms. 
The third inverse problem IP3 will be treated in Section \ref{th3} and concluding remarks will be made in the final Section \ref{th4}.

%Finally, the increasing stability of time-dependent source term can be obtained by using the boundary Dirichlet data.

\section{Main results}\label{man}
In this paper we always assume that the source function $F$ is required to be real-valued and the Fourier transform of $F$ is defined as
\[\hat{F}(\xi,\omega):={(2\pi)}^{-2}\int_{\R^4} F(x,t)e^{-i(\xi\cdot x+\omega t)}dxdt,\]
which implies that $\hat{F}(-\xi,-\omega)=\overline{F}(\xi,\omega)$ for all $\xi\in\R^3$ and $\omega\in \R$.

First we show increasing stability for the time-dependent inverse problem (IP1). 
Let $F(x,t)=f(x)g(t)$. It is supposed that $f(x)\in H^1(\R^3)$ and $g(t)\in H^2(0,\infty)$, where $f$ is compactly supported in $B_R$ and $g$ is supported in $(0,T_0)$ for some $T_0>0$.
The one-dimensional Fourier transform of $g$ with respect to the time variable $t$ is defined as follows:
\[\hat{g}(\omega)={(2\pi)}^{-\frac{1}{2}}\int_{\R} g(t)e^{-i\omega t}dt={(2\pi)}^{-\frac{1}{2}}\int_0^{T_0} g(t)e^{-i\omega t}dt.\]
We suppose there exist a number $b>1$ and a constant $\delta>0$ such that 
\be\label{af}
|\hat{g}(\omega)|\geq \delta>0 \quad\mbox{for all}\quad \omega\in (0,b).
\en
Physically, the parameter $b$ in \eqref{af} is associated with the bandwidth of the temporal signal $g(t)$. The condition \eqref{af} covers a large class of functions. For example, if 
$g(t)=e^{-(t-1)^2/\eta}\chi(t)$ with some $\eta>0$ and $\chi(t)\in C_0^\infty(\R)$ such that $\chi(t)=0$ for  $t\notin[0, T_0]$, the one can always find the parameters $b>0$ and $\delta>0$ such that \eqref{af} holds true. 
Since the source function is real valued, \eqref{af} holds true for all $\omega\in(-b, b)$.
We remark that the interval $(0, b)$ in \eqref{af} can be replaced by $(\omega_0-b, \omega_0+b)$ for some $\omega_0\geq 0$. In this paper we take $\omega_0=0$ for simplicity.

Through this paper $C>0$ denotes a generic constant which is independent of $b$, $\epsilon$ and $\Lambda$, may vary from line to line.
In the following theorem, we establish the increasing stability estimate of the $L^2$-norm of $f$ in terms of the parameter $b$ specified in (\ref{af}).  
\begin{theorem}\label{TH1} 
Let the condition \eqref{af} hold and
let $T>2R+T_0$.  Assume that $g(t)$ is given and  $\|f\|_{H^1(\R^3)}\leq M$ where $M>1$ is a constant.  Then
\begin{equation}\label{es0}
   \|f\|_{L^2(\R^3)}^2\leq  C(b^5\epsilon^2+\frac{M^2}{b^{\frac{4}{3}}|\ln \epsilon|^{\frac{1}{2}}}),
\end{equation}
where $\epsilon=\|u\|_{H^2\big([0,T];H^\frac{3}{2}(\partial B_R)\big)}$.
\end{theorem}

\begin{remark}
There are two parts in the stability estimate (\ref{es0}): the first  parts are the data discrepancy, while the second part comes from the high frequency tail of the function. 
The coefficient appearing in the Lipschitz part of $(\ref{es0})$ is polynomial type. However,  since $b$ is fixed in practice,  these coefficients are constants and do not pose any problem. It is clear to conclude that the ill-posedness of the inverse time-dependent source problem decreases as the parameter $b$ increases. The estimate (\ref{es0}) is in consistent with the increasing stability results of \cite{Li2017Increasing} in the frequency domain: the ill-posedness decreases when the width of the wave number interval $(0, b)$ increases. 
\end{remark} 
Next, we consider the second inverse problem (IP2).
Our aim is to establish an increasing stability estimate of $L^2$-norm of $F$ in terms of $\Lambda$.
\begin{theorem}\label{TH3} 
Let $T>\frac{2R}{\lambda}+T_0$ and let $F(x,t)\in H^3([0,T_0]; H^1(\R^3))$ be such that supp $F(x,t)$ $\subset B_R\times(0,T_0)$. Assume $\|F\|_{H^3([0,T_0]; H^1(\R^3))}\leq M$ for some $M>1$.
Then there exist constants $C>0$ and $\alpha\in(0,1)$ such that
 \begin{align}\label{es2}
 \|F\|^2_{L^2(\R^4)} \leq  C\Big(\Lambda^{10}\,\epsilon^2+\frac{M^2}{\Lambda|\ln\epsilon|^{\frac{1}{2}(1-\alpha)}}\Big),
 \end{align}
where $\epsilon=\sup\limits_{0<\lambda<\Lambda^2}\|u_\lambda\|_{H^2([0,T];H^\frac{3}{2}(\partial B_R))}$.
\end{theorem}
Finally, we present the increasing stability for the third inverse problem (IP3).
%Let $\lambda=1$ and assume that $F(x,t)$ takes the form
%\begin{equation}\label{3F}
%F(x,t)=f(\tilde{x},t)g(x_3), \quad \tilde{x}\in\R^2,  \quad   %t\in(0,+\infty),
%\end{equation}
%where $f$ is compactly supported on $\widetilde{B}_{R_0}\times(0,T_0)$ for some $0<R_0<R/\sqrt{2}$ and $g\in H^1(\R)$ is supported on $(-R_0,R_0)$. 
%where $\tilde{x}=(x_1,x_2)\in\R^2$ for $x=(x_1,x_2,x_3)\in\R^3$ and $\widetilde{B}_{R_0}:=\{\tilde{x}\in\R^2|\ |\tilde{x}|<R_0\}$. 
Suppose that $g$ is given and supported in $(-R_0, R_0)$ for some $0<R_0<R/\sqrt{2}$. 
The Fourier transform of  $g(x_3)$ is given by
\begin{equation}\nonumber
\hat{g}(\xi_3)=(2\pi)^{-\frac{1}{2}}\int_{\R}g(x_3)e^{-i\xi_3\cdot x_3}dx_3.
\end{equation}
We suppose
\be\label{g}
|\hat{g}(\xi_3)|\geq \delta>0 \quad\mbox{for all}\quad \xi_3\in (-b,b),
\en
where $b>1$.
\begin{theorem}\label{TH4}
%Suppose that the boundary surface measurement data is positive in the sense that
%\be
%\epsilon=||u||_{H^2([0,T];H^\frac{3}{2}(\partial B_R))}\geq \epsilon_0>0,
%\en 
%where $\epsilon_0>0$ is given.
Let $T>2R+T_0$ and let $f(\tilde{x},t)\in H^2([0,T_0]; H^1(\R^2))$ be such that supp$f(\tilde{x},t)\subset \widetilde{B}_{R_0}\times(0,T_0)$.
Assume $\|f\|_{H^1(\widetilde{B}_{R_0}\times(0,T_0))}\leq M$ for some $M>1$. Then there exist constants $C>0$ and $\alpha\in(0,1)$ such that
\begin{align}\label{es3}
\|f\|^2_{L^2(\R^3)} \leq C\Big(b^5\epsilon^2+b^3\,\,e^{2b(1-\alpha)}(1+b)^{2\alpha}\epsilon^{2\alpha}
+\frac{M^2}{b^{\frac{4}{3}}|\alpha\ln\epsilon|^{\frac{1}{2}}}\Big),
\end{align}
where $\epsilon=\|u\|_{H^2([0,T];H^\frac{3}{2}(\partial B_R))}$.
\end{theorem}

\section{Preliminaries}\label{pre}
Although partial results of this section are well known, we still present them for the readers' convenience and also to make this paper self-contained. Consider the time-dependent wave propagation problem caused by the source term $f(x)g(t)$. The mathematical model is
\begin{align}\label{eq2.1}
 \begin{cases}
 \partial_t^2u-\Delta u =f(x)g(t),\qquad  &(x,t)\in  \mathbb{R}^3\times(0,+\infty),\\
u(x,0)=\partial_tu(x,0)=0,\qquad  &x\in\mathbb{R}^3.
\end{cases}
\end{align}
%Specifically, the source function $F$ is assumed to be given in the following form
%\[F(x,t)=f(t)g(x),\]
If $f(x)\in L^2(\R^3)$ is  compactly supported in $B_R$ and 
$g(t)\in L^2(0,\infty)$ is supported in $(0, T)$. Then the problem (\ref{eq2.1}) admits a unique solution
\begin{equation}
 \nonumber
  u(x,t)\in{\mathcal{C}([0,+\infty);H^1(\mathbb{R}^3))\cap \mathcal{C}^1([0,+\infty);L^2(\mathbb{R}^3))}.
 \end{equation}
It's easy to prove this result by using the elliptic regularity properties of the Laplace operator (see \cite{Hsiao2008Boundary,McLean2000Strongly}).
Now we state the regularity of the solution for the initial boundary value problem (\ref{eq2.1}).
Below we suppose that $T>0$ is sufficiently large.
\begin{lemma}\label{lem2}
For $p>0$, let $f(x)\in H^p(\mathbb{R}^3)$ be supported in $B_R$ and $g(t)\in L^2([0,T])$. Then the problem $(\ref{eq2.1})$ admits a unique solution $u(x,t)\in{\mathcal{C}([0,T];H^{p+1}(\mathbb{R}^3))\cap H^2([0,T];H^{p-1}(\mathbb{R}^3))}$ satisfying
\begin{equation}\label{eq2.4}
\|u\|_{\mathcal{C}([0,T];H^{p+1}(\mathbb{R}^3))}+\|u\|_{H^2([0,T];H^{p-1}(\mathbb{R}^3))}
\leqslant C\|g\|_{L^2[0,T]}\,\|f\|_{H^{p}(\mathbb{R}^3)}
\end{equation}
with positive constant C depending on $R$ and $T$.
\end{lemma}

\begin{proof}
Denote the Fourier transform of $u(x,t)$ with respect to spatial variables as following:
\[ \hat{u}(\xi,t)=(2\pi)^{-3/2}\int_{\R^3}u(x,t)e^{-i\xi \cdot x}dx,\quad \xi\in \R^3.\]
Since $\widehat{\Delta u}(x,t)=-|\xi|^2\hat{u}(\xi,t)$ for all $\xi\in\R^3$, the function $\hat{u}(\xi,t)$ solves

\begin{align}\label{eq2.5}
\begin{cases}
 \partial_t^2 \hat{u}(\xi,t) + |\xi|^2\hat{u}(\xi,t)=\hat{f}(\xi)g(t), \qquad   &(\xi,t)\in\R^3\times(0,T), \\
 \hat{u}(\xi,0)=0,\quad \partial_t\hat{u}(\xi,0)=0, \qquad  &\xi\in\R^3.
\end{cases}
\end{align}
By Duhamel's principle, it is easy to check that the unique solution to (\ref{eq2.5}) takes the form
\begin{equation}
 \nonumber
\hat{u}(\xi,t)=\hat{f}(\xi)\int_{0}^{t}|\xi|^{-1}\sin(|\xi|(t-s))g(s)\,
   \mathrm{d}s.
 \end{equation}
For all $t\in[0,T]$ and $s\in[0,t]$, we introduce a function
\begin{equation}
\nonumber
 H(\cdot,t-s):= \xi \mapsto |\xi|^{-1}\sin(|\xi|(t-s))\;\widehat{f}(\xi).
\end{equation}
Then we have
\begin{equation}
\nonumber
\hat{u}(\xi,t)=\int_{0}^{t}H(\xi ,t-s)g(s)\,
   \mathrm{d}s.
\end{equation}
Since $f(x)$ is compactly supported in $B_R$, we obtain
 \begin{equation}\label{eq2.6}
 \begin{split}
   \|H(\cdot,t-s)\|_{L^2(\mathbb{R}^3)}^2 & =\int_{\mathbb{R}^3}|\xi|^{-2}\sin^2(|\xi|(t-s))|\widehat{f}(\xi)|^2\,
   \mathrm{d}\xi \\
   &\leqslant \int_{B_1}|\xi|^{-2}\sin^2(|\xi|(t-s))|\widehat{f}(\xi)|^2\,
   \mathrm{d}\xi+\int_{\mathbb{R}^3 \backslash B_1}|\xi|^{-2}|\widehat{f}(\xi)|^2\, \mathrm{d}\xi\\
   &\leqslant
   C|B_1|^2\|f\|_{L^2(\mathbb{R}^3)}^2+||f||_{L^2(\mathbb{R}^3)}^2\\
   &\leqslant(1+C|B_1|^2)\|f\|_{L^2(\mathbb{R}^3)}^2
 \end{split}
 \end{equation}
 for some positive constant $C$ and $B_1:=\{x\in\R^3|\ |x|<1\}$. In the same way, we know for $p>0$ that
 \begin{equation}\label{eq2.7}
 \begin{split}
  \|(1+|\xi|^2)^{\frac{p+1}{2}}H(\cdot,t-s)\|_{L^2(\mathbb{R}^3)}^2&=\int_{\mathbb{R}^3}(1+|\xi|^2)^{p+1}|\xi|^{-2}\sin^2(|\xi|(t-s))|\widehat{f}(\xi)|^2\,
   \mathrm{d}\xi \\
  &\leqslant \int_{B_1}(1+|\xi|^2)^{p+1}|\xi|^{-2}\sin^2(|\xi|(t-s))|\widehat{f}(\xi)|^2\,
   \mathrm{d}\xi\\
   &\qquad\qquad+\int_{\mathbb{R}^3 \backslash B_1}(1+|\xi|^2)^{p+1}|\xi|^{-2}|\widehat{f}(\xi)|^2\, \mathrm{d}\xi\\
    &\leqslant C\int_{\mathbb{R}^3}(1+|\xi|^{2})^p|\widehat{f}(\xi)|^2\,
   \mathrm{d}\xi+2\int_{\mathbb{R}^3 \backslash B_1}(1+|\xi|^2)^{p+1}\frac{1}{2|\xi|^2}|\widehat{f}(\xi)|^2\, \mathrm{d}\xi\\
    &\leqslant C\int_{\mathbb{R}^3}(1+|\xi|^{2})^p|\widehat{f}(\xi)|^2\,
   \mathrm{d}\xi+2\int_{\mathbb{R}^3}(1+|\xi|^{2})^p|\widehat{f}(\xi)|^2\,
   \mathrm{d}\xi\\
   &\leqslant (C+2)\|f\|_{H^p(\mathbb{R}^3)}^2.
   \end{split}
 \end{equation}
In view of estimates (\ref{eq2.6}) and (\ref{eq2.7}), it's easy to deduce that $u(x,t)\in \mathcal{C}([0,T];H^{p+1}(\mathbb{R}^3))$. Next we consider its
 derivative of time. Since for almost every $\xi\in\mathbb{R}^3$ we have
\begin{equation}
 \nonumber
 \partial_t \hat{u}(\xi,t)=\int_{0}^{t}\partial_t H(\xi,t-s)g(s)\, \mathrm{d}s=\int_{0}^{t}\cos(|\xi|(t-s))\;\hat{f}(\xi)g(s)\, \mathrm{d}s
\end{equation} 
and 
\begin{equation}
\nonumber
  \partial_t^2 \hat{u}(\xi,t)=g(t)\hat{f}(\xi)+\int_{0}^{t}\partial_t^2 H(\xi,t-s)g(s)\,\mathrm{d}s
  =g(t)\hat{f}(\xi)+\int_{0}^{t}-|\xi|\sin(|\xi|(t-s))\;\hat{f}(\xi)g(s)\,\mathrm{d}s.
\end{equation}
Combining this with (\ref{eq2.7}), we find 
\begin{equation}\label{eq2.8}
 \nonumber
 \|(1+|\xi|^2)^{\frac{p}{2}}\partial_t H(\xi,t-s)\|_{L^2(\mathbb{R}^3)}^2\leqslant  C\|f(x)\|_{H^p(\mathbb{R}^3)}^2 
\end{equation}
and
\begin{equation}
 \nonumber 
  \|(1+|\xi|^2)^{\frac{p-1}{2}}\partial_t^2 H(\xi,t-s)\|_{L^2(\mathbb{R}^3)}^2\leqslant  C\|f(x)\|_{H^p(\mathbb{R}^3)}^2.
 \end{equation}
From this fact the result $u(x,t)\in H^2([0,T];H^{p-1}(\mathbb{R}^3))$ follows. The estimate (\ref{eq2.4}) is also easy to get from the previous estimates.
\end{proof}

In Theorem \ref{TH1}-\ref{TH4}, we obtain increasing stability results by using only the Dirichlet boundary data rather than the Cauchy data.  This is due to the fact that our inverse problems are all formulated in the unbounded domain $|x|>R$ and the Neumann data on $\partial B_R\times(0,T)$ can be controlled by the Dirichlet data on $\partial B_R\times(0,T)$. To prove this rigorously, we need to impose higher regularity assumptions on $f$ and $g$.
\begin{lemma}\label{lem4}
Let $f(x)\in H^1(\mathbb{R}^3)$ and $g(t)\in H^2([0,T])$. Then the problem $(\ref{eq2.1})$ admits a unique solution $u(x,t)\in{\mathcal{C}^2([0,T];H^2(\mathbb{R}^3))\cap H^4([0,T];L^2(\mathbb{R}^3))}$ satisfying the estimate
\begin{equation}\label{eq2.9}
 \|\partial_{\nu}u\|_{L^2([0,T];H^{\frac{1}{2}}(\partial B_R))}
\leqslant C\|u\|_{H^2([0,T];H^{\frac{3}{2}}(\partial B_R))}
\end{equation}
with positive constant C depending on $T$ and $R$.
\end{lemma}

\begin{proof}
Applying Lemma \ref{lem2} with $p=1$, we know 
\be\nonumber
&u(x,t)\in{\mathcal{C}^2([0,T];H^2(\mathbb{R}^3))\cap H^4([0,T];L^2(\mathbb{R}^3))}, \\\label{h} &h(x,t):=u|_{\partial B_R\times[0,T]}\in\mathcal{C}^2([0,T];H^{\frac{3}{2}}(\partial B_R)).
\en 
Therefore, the restriction of $u(x,t)$ to $({\mathbb{R}^3\backslash\overline{B_R}})\times [0,T]$ solves the initial boundary value problem
\begin{align}\label{eq2.10}
\begin{cases}
  \partial_t^2 u(x,t) - \Delta u(x,t)=0,  \qquad   &(x,t)\in{\mathbb{R}^3\backslash\overline{B_R}}\times(0,T),\\
u(x,0)=\partial_tu(x,0)=0, \qquad    &x\in\mathbb{R}^3\backslash\overline{B_R}, \\
u(x,t)=h(x,t), \qquad    &(x,t)\in \partial B_R \times (0,T).
\end{cases}
\end{align}
Combining a classical lifting result with the fact that $h(x,0)=h_t(x,0)=0$, we introduce a function $H(x,t)\in \mathcal{C}^2([0,T];H^2({\mathbb{R}^3\backslash B_R}))$ such that $H|_{\partial B_R\times[0,T]}=h$, $H(x,0)=\partial_tH(x,0)=0$ and
\begin{equation}\label{eq2.11}
\|H\|_{H^2([0,T];H^2({\mathbb{R}^3\backslash B_R}))}\leqslant C\|h\|_{H^2([0,T];H^{\frac{3}{2}}(\partial B_R))},
\end{equation}
where $C>0$ depends on $R$ and $T$.
Therefore, we can split $u$ to $u=H+V$ on ${\mathbb{R}^3\backslash B_R}\times(0,T)$. Here $V$ solves
\begin{align*}
\begin{cases}
 \partial_t^2 V(x,t) - \Delta V(x,t)=-(\partial_t^2 H(x,t) - \Delta H(x,t)):=G, \quad &(x,t)\in{\mathbb{R}^3\backslash\overline{B_R}}\times(0,T), \\
V(x,0)=V_t(x,0)=0, \quad &x\in\mathbb{R}^3\backslash\overline{B_R}, \\
V(x,t)=0,\quad &(x,t)\in\partial B_R \times (0,T).
\end{cases}
\end{align*}
Using the fact that $G\in L^2([0,T];L^2({\mathbb{R}^3\backslash B_R})$ and $H(x,0)=0$ with Lemma \ref{lem2}, we get $V\in{\mathcal{C}([0,T];H^1(\mathbb{R}^3))\cap H^2([0,T];H^{-1}(\mathbb{R}^3))}$ satisfying the estimate
\begin{equation}
\nonumber
\|V\|_{H^2([0,T];H^{-1}(\mathbb{R}^3))}\leq C\|G\|_{L^2([0,T];L^2({\mathbb{R}^3\backslash \overline{B_R}})}\leq C\|H\|_{H^2([0,T];H^2({\mathbb{R}^3\backslash \overline{B_R}})}.
\end{equation}
Combining this with (\ref{eq2.11}), we deduce that
\begin{equation}\nonumber
\|u\|_{L^2([0,T];H^2({\mathbb{R}^3\backslash \overline{B_R}})}\leqslant C\|h\|_{H^2([0,T];H^{\frac{3}{2}}(\partial B_R))}
\end{equation}
and using the continuity of the trace map, we obtain
\begin{equation}\nonumber
\|\partial_{\nu}u\|_{L^2([0,T];H^{\frac{1}{2}}(\partial B_R))}
\leqslant C\|u\|_{L^2([0,T];H^2({\mathbb{R}^3\backslash\overline{B_R}}))}.
\end{equation}
Combining the last two estimates with the definition of $h$ given by \eqref{h}, we finally obtain (\ref{eq2.9}).
\end{proof}
For $r>0$, we denote $B(0,r)=\{x\in\R^d|\ |x|<r\}$. Below we state a stability estimate for analytic continuation problems, which can be seen in \cite{Mourad-Ben, Vessella}.
\begin{proposition}\label{pro3.1}
Let $\mathcal{O}$ be a non empty open set of the unit ball $B(0,1)\subset
\mathbb{R}^{d}$, $d\geq2$, and let $G$ be an analytic function in $B(0,2),$ that satisfy
$$\|\partial^{\gamma}G\|_{L^{\infty}(B(0,2))}\leq M_0\,|\gamma|!\;{\eta^{-|\gamma|}},\,\,\,\,\forall\,\gamma\in(\mathbb{N}\cup\{0\})^{d},$$
for some  $M_0>0$ and $\eta>0$. Then, we have
$$\|G\|_{L^{\infty}(B(0,1))}\leq N\, M_0^{1-\mu}\;\|G\|_{L^{\infty}(\mathcal{O})}^{\mu},$$
where $\mu\in(0,1)$ depends on $d$, $\eta$ and $|\mathcal{O}|$ and $N=N(\eta)>0$.
\end{proposition}

\section{Proofs of Theorem \ref{TH1}}\label{th1}
%To derive the increasing stability estimate of the source term $g$, we need to impose more restrictive assumptions on the temporal function  $f(t)\in H^2[0,T_0]$. We suppose there exist  a number $b>1$ and a constant $\delta=\delta(b)>0$ such that 
In this section, we discuss increasing
stability of the inverse problem IP1.
Firstly we introduce the time-dependent test function
\begin{equation}\nonumber
w(x,t):= e^{-i\xi\cdot x-i\omega t},
\end{equation}
where $\xi\in\mathbb{R}^3$ and $\omega\in\mathbb{R}$ satisfy $|\xi|^2=\omega^2$. It is obvious that $w$ satisfies the homogeneous acoustic wave equation
\begin{equation}\label{eq2.12}
\partial_t^2 w -  \Delta w =0  \qquad \textrm{in} \  \mathbb{R}^3\times(0,+\infty).
\end{equation}
Then, multiplying $w$ on both sides of the equation (\ref{eq1}) and integrating over $B_R\times(0,T)$, we obtain
\begin{equation}\label{eq2.13}
  \int_0^T\int_{B_R} (\partial_t^2u-\Delta u) w\ \mathrm{d}x\mathrm{d}t=\int_0^T\int_{B_R} f(x) g(t) w(x,t)\ \mathrm{d}x\mathrm{d}t.
\end{equation}
Using (\ref{eq2.12}) and integrating by parts, one deduces from the left hand side of (\ref{eq2.13}) that
\begin{align*}
  \int_0^T\int_{B_R}&(\partial_{t}^2u(x,t)-\Delta u (x,t)) w(x,t)\ \mathrm{d}x\mathrm{d}t \\
  =&\int_0^T\int_{B_R} \big(\partial_{t}^2u(x,t) w(x,t)- u (x,t) \partial_t^2w(x,t)\big) \ \mathrm{d}x\mathrm{d}t\\
  &\quad+\int_0^T\int_{\partial B_R} \big(u(x,t) \frac{\partial w(x,t)}{\partial \nu}-w(x,t) \frac{\partial u(x,t)}{\partial \nu}\big)\ \mathrm{d}s(x)\mathrm{d}t\\
  =&\int_{B_R} \big(\partial_{t}u(x,t) w(x,t)- u (x,t) \partial_{t}w(x,t)\big)\big|_0^T\;\mathrm{d}x\\
  &\quad+\int_0^T\int_{\partial B_R} \big(u(x,t) \frac{\partial w(x,t)}{\partial \nu}-w(x,t) \frac{\partial u(x,t)}{\partial \nu}\big)\ \mathrm{d}s(x)\mathrm{d}t.\\
=&\int_0^T\int_{\partial B_R} \big(u(x,t) \frac{\partial w(x,t)}{\partial \nu}-w(x,t) \frac{\partial u(x,t)}{\partial \nu}\big)\ \mathrm{d}s(x)\mathrm{d}t.
\end{align*}
Note that, in the last step we have used the fact that $u(x,t)=0$ when $|x|<R$ and $t>T_0+2R$, which follows straightforwardly from Huygens' principle (see 
 \cite[Lemma 2.1.]{BaoHu}). This implies $u(x,T)=\partial_t u(x,T)=0$ for $x\in B_R$ and $T>T_0+2R$. Hence, the integral over $B_R$ on the left hand side of the previous identity vanishes. Recalling the estimate in Lemma \ref{lem4} with Sobolev's embedding theorems, we bound the left hand side of (\ref{eq2.13}) by
\begin{align*}
  &\Big|\int_0^T\int_{B_R} (\partial_{t}^2u(x,t)-\Delta u (x,t)) w(x,t)\ \mathrm{d}x\mathrm{d}t \Big|\\
  = &\Big|\int_0^T\int_{\partial B_R} (u(x,t) \frac{\partial w(x,t)}{\partial \nu}-w(x,t) \frac{\partial u(x,t)}{\partial \nu})\ \mathrm{d}s(x)\mathrm{d}t\Big|\\
 \leq &\|u\|_{L^2([0,T];L^2(\partial B_R))}\|\frac{\partial w}{\partial \nu}\|_{L^2([0,T];L^2(\partial B_R))}
  +\|\frac{\partial u}{\partial \nu}\|_{L^2([0,T];L^2(\partial B_R))}\|w\|_{L^2([0,T];L^2(\partial B_R))}\\
  \leq &C(\|u\|_{L^2([0,T];L^2(\partial B_R))}\|w\|_{L^2([0,T];H^2(B_R))}
  +\|u\|_{H^2([0,T];H^\frac{3}{2}(\partial B_R))}\|w\|_{L^2([0,T];L^2(\partial B_R))})\\
  \leq &C\|u\|_{H^2([0,T];H^\frac{3}{2}(\partial B_R))} \|w\|_{L^2([0,T];H^2(B_R))}.\\
\end{align*}
By the definition of $w$, one can check that
\ben
\|w\|_{L^2([0,T];H^2(B_R))}\leq C\, (1+|\xi|)\qquad\mbox{for all}\quad |\xi|^2=\omega^2.
\enn
Hence, using the assumption about $g$ and the fact that $f$ is supported in $B_R$ , we derive 
from (\ref{eq2.13}) together with the previous two relations that
\ben
|\hat{f}(\xi)\hat{g}(\omega)|
&=&\Big|(2\pi)^{-2}\int_0^T\int_{\mathbb{R}^3} f(x)g(t)  e^{-i\xi\cdot x-i\omega t}\ \mathrm{d}x\mathrm{d}t\Big|\\
 &=&\Big|(2\pi)^{-2}\int_0^T\int_{B_R} f(x)g(t) w(x,t)\ \mathrm{d}x\mathrm{d}t\Big|\\
&\leq &C(1+|\xi|)\|u\|_{H^2([0,T];H^\frac{3}{2}(\partial B_R))}.
\enn
In view of the assumption \eqref{af}, one obtains for $\omega\in(0, b)$ and  $|\xi|^2=\omega^2$ that
\begin{equation}\label{eq2.14}
 |\hat{f}(\xi)|\leq C\frac{(1+|\xi|^2)\|u\|_{H^2([0,T],H^\frac{3}{2}(\partial B_R))}}{|\hat{g}(\omega)|}\leq C\,\delta^{-1} (1+|\xi|) \|u\|_{H^2([0,T];H^\frac{3}{2}(\partial B_R))}.
\end{equation}
We note that (\ref{eq2.14}) gives an estimate of $\hat{f}(\xi)$ over the domain $E:=B(0, b)=\{\xi\in\mathbb{R}^3| \  |\xi|<b\}$, that is,
\be\label{2.14}
\|\hat{f}\|_{L^\infty(E)}\leq C\,\delta^{-1} (1+b)\, \epsilon,
\en where $\epsilon=\|u\|_{H^2([0,T];H^\frac{3}{2}(\partial B_R))}$ represents the measurement data on $\partial B_R\times (0, T)$.

%\subsection{Proof of Theorem \ref{TH1}}%

For $f(x)\in L^2(\mathbb{R}^3)$ and  $\mathrm{supp}f\subset B_R $. Applying the Parseval's identity, we have 
 \begin{equation}
 \nonumber
  \|f\|_{L^2(\mathbb{R}^3)}^2=\|\hat{f}\|_{L^2(\mathbb{R}^3)}^2=I(k)+\int_{|\xi|>k}|\hat{f}(\xi)|^2\,\mathrm{d}\xi,
 \end{equation}
where
\begin{align}\label{I}
   I(k):=   \int_{|\xi|\leq k} |\hat{f}(\xi)|^2\,\mathrm{d}\xi&=\int_{0}^{k}\int_{\mathbb{S}^2}|\hat{f}(l\theta)|^2 l^2 \mathrm{d}\theta\mathrm{d}l.
\end{align}
Obviously, the following inequalities hold
\begin{align*}
\int_{|\xi|>k}|\hat{f}(\xi)|^2\,\mathrm{d}\xi\leq \frac{1}{k^2}\,\int_{|\xi|>k}|\widehat{\nabla f}(\xi)|^2\,\mathrm{d}\xi\leq\frac{M^2}{k^2}
\end{align*}
by the Parseval's identity.
Since the integrand is an entire analytic function of $\xi$, the integral $I(k)$ with respect to $\xi$ can be taken over by any path joining points $0$ and $k$ in complex plane. Thus $I(k)$ is an entire analytic function of $k=k_1+ik_2$ $(k_1,k_2\in \mathbb{R})$ and the following estimate holds.
\begin{lemma}\label{1}
 Let $f(x)\in L^2(\mathbb{R}^3)$ and $\mathrm{supp}\ f\subset B_R $. Then
 \begin{equation}\label{esI}
   |I(k)|\leq (\frac{4\pi}{3})^2\,R^3\,|k|^3e^{2R|k_2|}\|f\|^2_{L^2(\mathbb{R}^3)},\qquad k=k_1+i k_2\in \C.
 \end{equation}
\end{lemma}
\begin{proof}
Set $l=ks$ for $s\in (0,1)$. Then it is easy to get
\begin{align*}
      I(k)&=\int_{0}^{k}\int_{\mathbb{S}^2} |\hat{f}(l\theta)|^2 l^2 \mathrm{d}\theta\mathrm{d} l\\
      &=\int_{0}^{1}\int_{\mathbb{S}^2}|\hat{f}(ks\theta)|^2 k^3s^2 \mathrm{d}\theta\mathrm{d}s.
\end{align*}
Noting the elementary inequality $|e^{2k_2s\theta\cdot x}|\leq e^{2R|k_2|}$ for all $x\in B_R$ and $\theta\in \mathbb{S}^2$, we have
\begin{align*}
|I(k)|&=\left|\int_{0}^{1}\int_{\mathbb{S}^2} \left|\int_{B_R}f(x)e^{-iks\theta\cdot x}\mathrm{d} x\right|^2 k^3s^2 \mathrm{d}\theta\mathrm{d} s\right|\\
&\leq \frac{4\pi}{3}R^3\int_{0}^{1}\int_{\mathbb{S}^2}|k|^3s^2 \left(\int_{B_R}|f(x)|^2|e^{2k_2s\theta\cdot x}|\, \mathrm{d} x\right)\mathrm{d}\theta\mathrm{d} s\\
&\leq (\frac{4\pi}{3})^2\,R^3\,|k|^3e^{2R|k_2|}\|f\|^2_{L^2(\mathbb{R}^3)}.
\end{align*}
This completes the Lemma \ref{1}.
\end{proof}

The following Lemma is essential to show the relation between $I(k)$ for $k\in(b,\infty)$ with $I(b).$ Its proof can be found in \cite{Cheng2016Increasing}.
\begin{lemma}\label{2}
  Let $J(z)$ be an analytic function in $S=\{z=x+iy\in \mathbb{C}:-\frac{\pi}{4}<\arg z<\frac{\pi}{4}\}$ and continuous in $\overline{S}$ satisfying
  \begin{equation}
  \nonumber
  \begin{cases}
    |J(z)| \leq \epsilon, \  & z\in (0,\ L], \\
    |J(z)| \leq V, \  & z\in S,\\
    |J(0)|  =0.
  \end{cases}
  \end{equation}
  Then there exists a function $\mu(z)$ satisfying
  \begin{equation}
  \nonumber
  \begin{cases}
   \mu (z)  \geq\frac{1}{2},\ \ &z\in (L,\ 2^{\frac{1}{4}}L), \\
   \mu (z)  \geq\frac{1}{\pi}((\frac{z}{L})^{4}-1)^{-\frac{1}{2}},\ \ & z\in (2^{\frac{1}{4}}L,\ \infty)
  \end{cases}
  \end{equation}
   such that
\begin{equation}
\nonumber
 |J(z)|\leq V\epsilon^{\mu(z)}, \ \ \forall z\in (L,\ \infty).
\end{equation}
\end{lemma}
Let the sector $S\subset \C$ be given as in Lemma  \ref{2}.
Now, it follows from Lemma \ref{1} that
\begin{equation}
\nonumber
|I(k)e^{-(2R+1)k}|\leq CM^2 \quad \mbox{for all}\quad k\in S.
\end{equation}
Recalling from a priori estimate \eqref{2.14}, we obtain 
\begin{equation}
\nonumber
I(k)=\int_{|\xi|\leq k}|\hat{f}(\xi)|^2d\xi\leq C\,k^3\,(1+k)^2\,\epsilon^2, \quad k\in(0, b].
\end{equation}
Hence
\begin{equation}
\nonumber
| I(k)e^{-(2R+1)k}|\leq C \epsilon^2, \ \ k\in(0,b].
\end{equation}
Then applying Lemma \ref{2} with $L=b$ to the function $J(k):=I(k)e^{-(2R+1)k}$, we know that there exists a function $\mu(k)$ satisfying
  \begin{equation}\
  \begin{cases}
  \mu (k)  \geq\frac{1}{2},\ \ &k\in (b,\ 2^{\frac{1}{4}}b), \\
   \mu (k)  \geq\frac{1}{\pi}((\frac{k}{b})^{4}-1)^{-\frac{1}{2}},\ \ &k\in (2^{\frac{1}{4}}b,\ \infty)
  \end{cases}
\end{equation}
   such that
  \begin{equation}
   \nonumber
   |I(k)e^{-(2R+1)k}|\leq CM^2\epsilon^{2\mu} \ \ \quad\mbox{for all}\quad k\in(b,\infty).
   \end{equation}
Now we show the proof of Theorem \ref{TH1}. We  assume that $\epsilon<e^{-1}$, otherwise the estimate is obvious. 
Let
\begin{equation}
\nonumber
k= \begin{cases}
  \frac{1}{((2R+3)\pi)^{\frac{1}{3}}}b^{\frac{2}{3}}|\ln \epsilon|^{\frac{1}{4}},  \qquad&\mbox{if} \quad 2^{\frac{1}{4}} ((2R+3)\pi)^{\frac{1}{3}}b^{\frac{1}{3}}< |\ln \epsilon|^{\frac{1}{4}},\\
     b,  \qquad&\mbox{if} \quad |\ln \epsilon|^{\frac{1}{4}}\leq 2^{\frac{1}{4}}((2R+3)\pi)^{\frac{1}{3}} b^{\frac{1}{3}}.
  \end{cases}
\end{equation}
Case (i):  $2^{\frac{1}{4}}((2R+3)\pi)^{\frac{1}{3}} b^{\frac{1}{3}}< |\ln \epsilon|^{\frac{1}{4}}$. Then we have
\begin{eqnarray*}
 %|I(k)|&\leq& CM^2\epsilon^{2\mu}e^{(2R+1)k} \\
 |I(k)|&\leq& CM^2\epsilon^{2\mu}e^{(2R+1)k}\\
 &=&CM^2e^{(2R+1)k-2\mu(k)|\ln\epsilon|}\\
           &\leq& CM^2e^{\frac{(2R+3)}{((2R+3)\pi)^{\frac{1}{3}}}b^{\frac{2}{3}}|\ln \epsilon|^{\frac{1}{4}}-\frac{2|\ln \epsilon|}{\pi}(\frac{b}{k})^2 }\\
           &=& CM^2 e^{-2\big(\frac{(2R+3)^2}{\pi}\big)^{\frac{1}{3}}b^{\frac{2}{3}}|\ln \epsilon|^{\frac{1}{2}}(1- \frac{1}{2}|\ln \epsilon|^{-\frac{1}{4}}) }. \\
\end{eqnarray*}
Noting that $|\ln\epsilon|^{-\frac{1}{4}}<1$ and $\big(\frac{(2R+3)^2}{\pi}\big)^{\frac{1}{3}}>1$, we have 
\[|I(k)|\leq CM^2e^{-b^{\frac{2}{3}}|\ln \epsilon|^{\frac{1}{2}}}.\]
Using the inequality $e^{-t}\leq\frac{6!}{t^6}$ for $t>0$, we get
\begin{equation}\label{ess}
  |I(k)|\leq CM^2\frac{1}{b^4|\ln \epsilon|^{3}}.
\end{equation}
Since $b^4|\ln \epsilon|^{3}\geq b^{\frac{4}{3}}|\ln \epsilon|^{\frac{1}{2}}$ when $b>1$ and $|\ln\epsilon|>1$.
Hence
\begin{equation}\label{ess1}
  \begin{split}
     \|f\|_{L^2(\mathbb{R}^3)}^2&= I(k)+\int_{|\xi|>k} |\widehat{f}(\xi)|^2\,\mathrm{d}\xi\\
     &\leq I(k)+\frac{M^2}{k^2}\\
                         &\leq C(\frac{ M^2}{b^4|\ln \epsilon|^{3}}+\frac{ M^2}{b^{\frac{4}{3}}|\ln \epsilon|^{\frac{1}{2}}})\\
                         &\leq \frac{ CM^2}{b^{\frac{4}{3}}|\ln \epsilon|^{\frac{1}{2}}}.
  \end{split}
\end{equation}
Case (ii):  $|\ln \epsilon|^{\frac{1}{4}}\leq 2^{\frac{1}{4}}((2R+3)\pi)^{\frac{1}{3}} b^{\frac{1}{3}}$. In this case we have $k=b$ by the choice of $k$ and $|I(b)|\leq C\,b^5\epsilon^2$. Hence,
 \begin{equation}\label{ess2}
  \begin{split}
     \|f\|_{L^2(\mathbb{R}^3)}^2&=I(b)+\int_{|\xi|>b} |\widehat{f}(\xi)|^2\,\mathrm{d}\xi\\
                         &\leq C(b^5\epsilon^2+\frac{ M^2}{b^2})\\
                         &\leq C(b^5\epsilon^2+\frac{M^2}{b^{\frac{4}{3}}|\ln \epsilon|^{\frac{1}{2}}}).
  \end{split}
 \end{equation}
Combining (\ref{ess1}) and (\ref{ess2}), we finally get
\begin{equation*}
\|f\|_{L^2(\mathbb{R}^3)}^2\leq  C(b^5\epsilon^2+\frac{M^2}{b^{\frac{4}{3}}|\ln \epsilon|^{\frac{1}{2}}}).
\end{equation*}

%\end{proof}

\section{Proof of Theorem \ref{TH3}}\label{th2}
\setcounter{equation}{0}

In this section, we consider the following initial value problem for the wave equation
\begin{equation}\nonumber
\begin{cases}
\partial_t^2u_\lambda(x,t)-\lambda\Delta u_\lambda(x,t)=F(x,t), \ &(x,t)\in\R^3\times(0,\infty),\\
u(x,0)=0, \quad \partial_tu(x,0)=0, \ &x\in\R^3.
\end{cases}
\end{equation}
Let $0<\lambda<\Lambda^2$. Our aim is to recover the compacted supported function $F(x,t)$ from the data $\{u_\lambda(x, t)|\ x\in\partial B_R, \ t\in(0,T)\}$. Physically, such kind of the measurement data
can be obtained by changing the background medium artificially and locally for the purpose of
recovering a time-dependent source term which might be non-radiating for a fixed parameter $\lambda$.
 
As done in last section we introduce the test function
 $$w(x,t)=e^{-i(\xi\cdot x+\omega t)},\qquad \omega^2-\lambda|\xi|^2=0.$$ One can check that $\partial_t^2w-\lambda\Delta w=0$. In this case we also have the strong Huygens' principle that $u_\lambda(x,t)=0$ when $|x|<R$ and $t>T_0+\frac{2R}{\lambda}$ (see \cite[Lemma 2.1.]{BaoHu}).

\begin{figure}[htbp]
 \centering	
\includegraphics[width=0.48\textwidth]{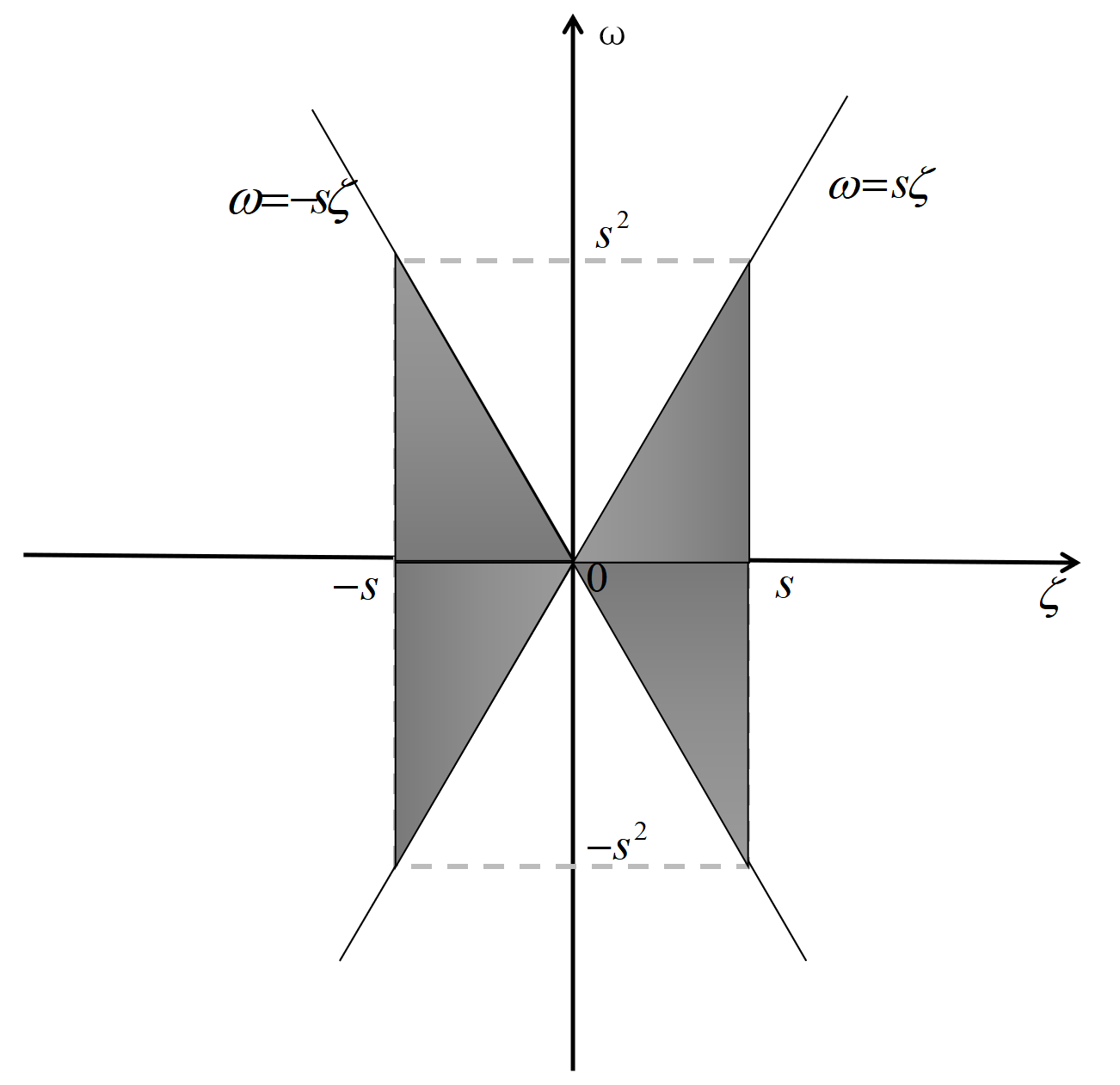}
	\caption{$E(s)$ is the shaded area, $|\xi|^2=\zeta^2$.}
	\label{Fig1}
\end{figure}

\begin{figure}[htbp]
 \centering	
\includegraphics[width=0.48\textwidth]{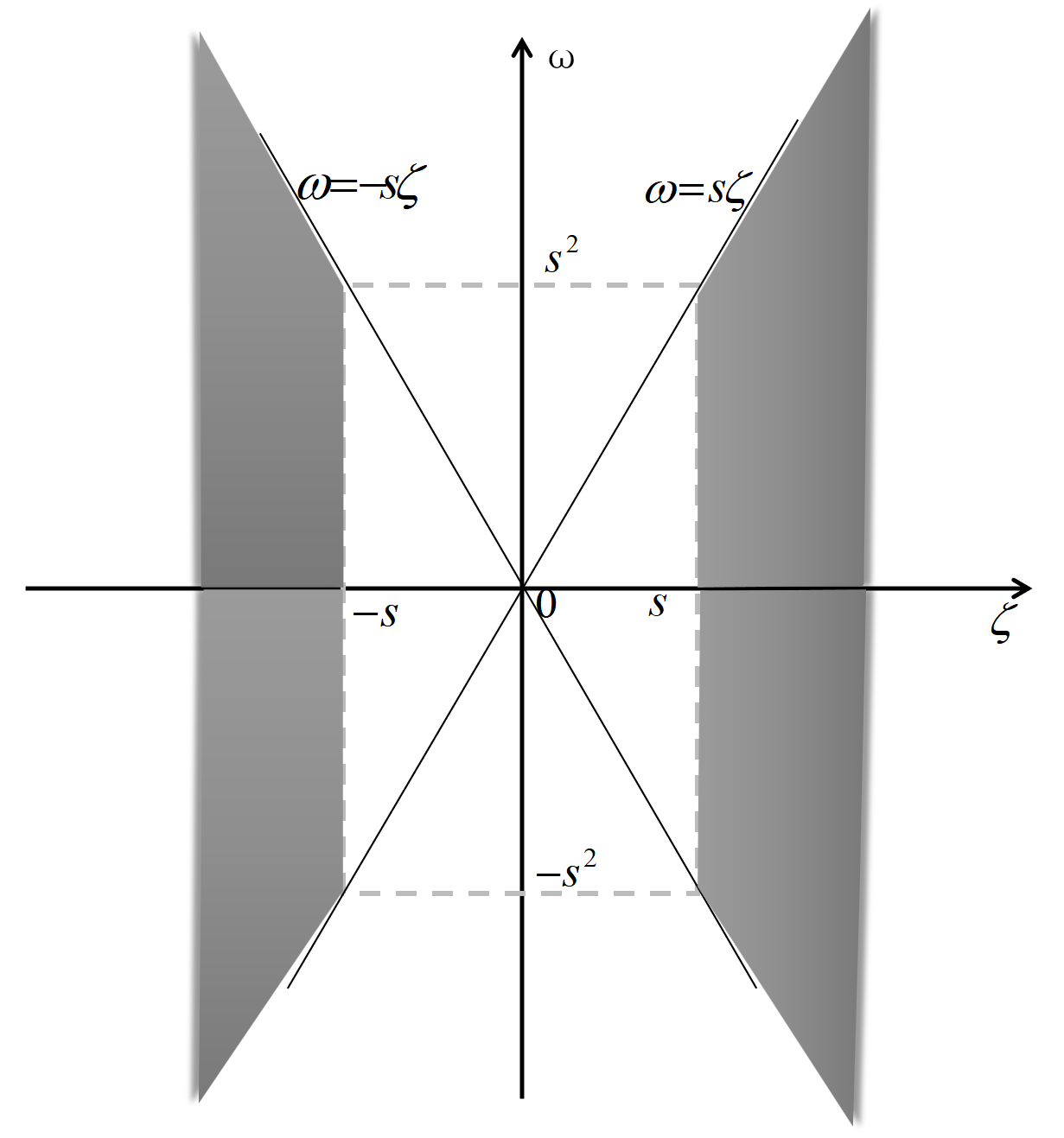}
	\caption{$E_1(s)$ is the shaded area, $|\xi|^2=\zeta^2$.}
	\label{Fig2}
\end{figure}

\begin{figure}[htbp]
 \centering	
\includegraphics[width=0.48\textwidth]{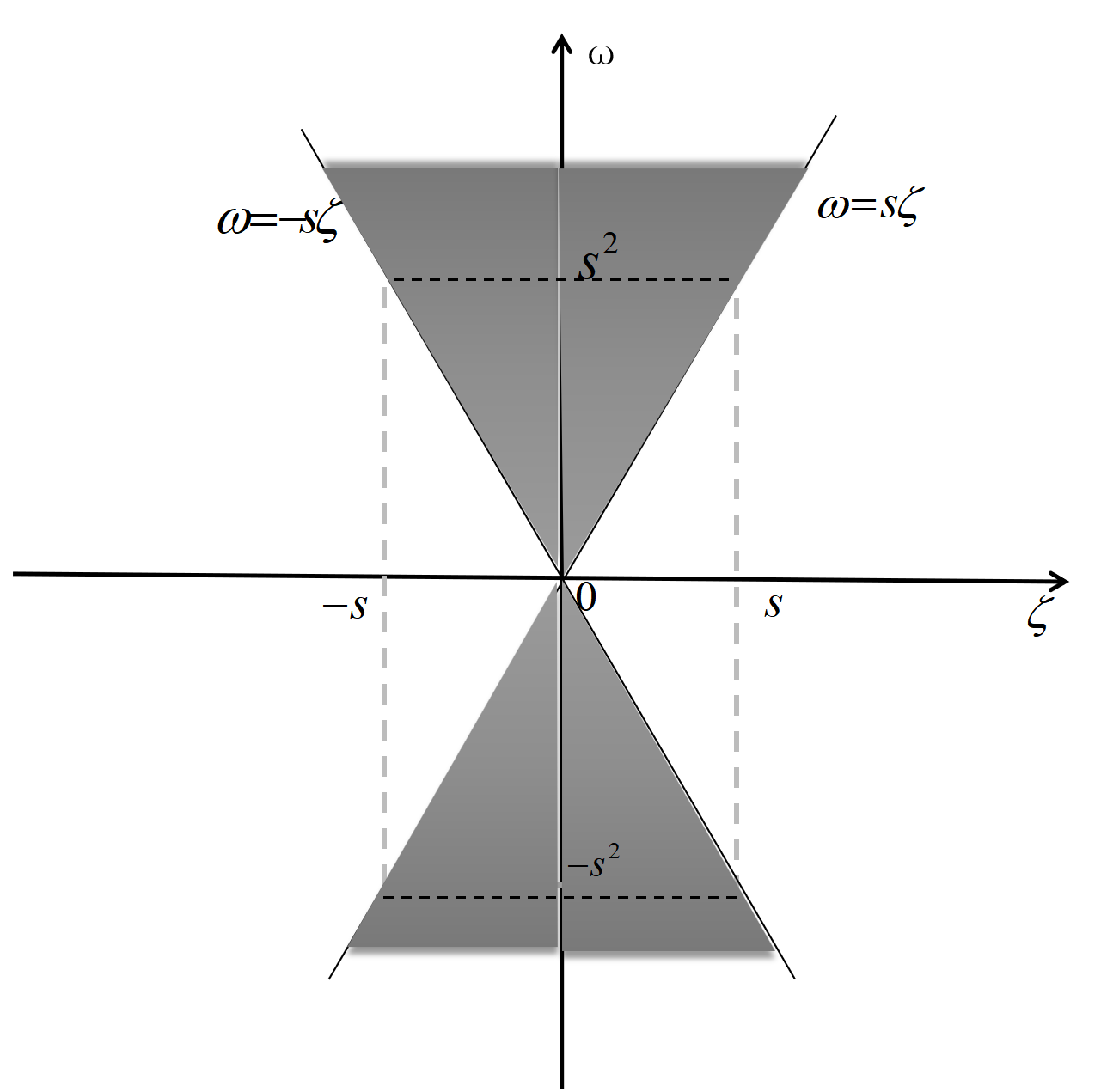}
	\caption{$E_2(s)$ is the shaded area, $|\xi|^2=\zeta^2$.}
	\label{Fig3}
\end{figure}

Multiplying $w$ on both sides of wave equation   (\ref{eq1}) and integrating over $B_R\times(0,T)$, we obtain
\begin{align*}
 \int_0^T\int_{B_R}F(x,t)e^{-i(\xi\cdot x+\omega t)}dxdt
 &=\int_0^T\int_{B_R}\big(\partial_t^2u_\lambda(x,t)-\lambda\Delta u_\lambda(x,t)\big) w(x,t)\ \mathrm{d}x\mathrm{d}t\\
 &=
  \lambda\int_0^T\int_{\partial B_R}\big(u_\lambda(x,t)\frac{\partial w(x,t)}{\partial \nu}-w(x,t) \frac{\partial u_\lambda(x,t)}{\partial \nu}\big)\ \mathrm{d}s(x)\mathrm{d}t.
\end{align*}
Define
\begin{equation}\nonumber
\hat{F}(\xi,\omega)=(2\pi)^{-2}\int_{\R^4}F(x,t)e^{-i(\xi\cdot x+\omega t)}dxdt,
\end{equation}
with $(\xi,\omega)=(\xi_1,\xi_2,\xi_3,\omega)\in\R^4$.
Since supp $F(x,t)$ $\subset B_R\times(0,T_0)$, we have
\be
\nonumber
(2\pi)^2\hat{F}(\xi,\omega)
&=&\int_0^T\int_{B_R}F(x,t)e^{-i(\xi\cdot x+\omega t)}dxdt \\ \label{Fu}
& =&\lambda\int_0^T\int_{\partial B_R}\big(u_\lambda(x,t)\frac{\partial w(x,t)}{\partial \nu}-w(x,t) \frac{\partial u_\lambda(x,t)}{\partial \nu}\big)\ \mathrm{d}s(x)\mathrm{d}t.
\en
Using the fact that $F(x,t)$ has compact support and arguing similarly to the proof of Lemma \ref{lem4},  we deduce
\begin{equation}\label{u}
\|\partial_{\nu}u_\lambda\|_{L^2([0,T];H^{\frac{1}{2}}(\partial B_{R}))}
\leqslant C\|u_\lambda\|_{H^2([0,T];H^{\frac{3}{2}}(\partial B_{R}))}.
\end{equation}
Consider the set (see Figure $\ref{Fig1}$)
\begin{equation}\nonumber
E(s)=\{(\xi,\omega)\in \R^4| \ \omega^2-\lambda|\xi|^2=0,\ |\xi|<s,\ 0<\lambda<s^2\}.
\end{equation}
Combining $(\ref{Fu})$ and $(\ref{u})$, we have
\begin{equation}\label{hatH}
|\hat{F}(\xi,\omega)|\leq C\lambda(1+|\xi|)\|u_\lambda\|_{H^2([0,T];H^\frac{3}{2}(\partial B_R))}\leq C\lambda(1+|\xi|)\epsilon, \quad (\xi,\omega)\in {E(\Lambda)},
\end{equation}
where $\epsilon=\sup\limits_{0<\lambda<\Lambda^2}\|u_\lambda\|_{H^2([0,T];H^\frac{3}{2}(\partial B_R))}$.
Denote
$I(s):=\int_{E(s)}|\hat{F}(\xi,\omega)|^2d\xi\,d\omega$,
then 
\begin{align}\label{IE}
|I(s)|\leq C\,|E(s)|s^4\,(1+s)^2\,\epsilon^2 \quad \mbox{for all}\quad s\in[0,\Lambda].
\end{align}
%Since 
%\begin{equation}
%\begin{split}
%\int_0^{s^2\hat{r}}|\hat{F}|^2d\omega
%=\int_0^1|\int_{\R^4}F(x,t)e^{i(s\hat{r}\hat{\xi}\cdot x+s^2\hat{r}\hat{\omega}t)}dxdt|^2s^2\hat{r}d\hat{\omega},
%\end{split}
%\end{equation}
Using the polar coordinates $\xi=r\hat{\xi}=r(\cos\theta\sin\varphi, \sin\theta\sin\varphi, \cos\varphi)$ for $0\leq\theta\leq2\pi$ and $0\leq\varphi\leq\pi$, we deduce that 
\begin{equation}\label{Ie}
\begin{split}
I(s)&=\int_{E(s)}|\hat{F}(\xi,\omega)|^2d\xi d\omega\\
&=\int_0^{2\pi}\int_0^{\pi}\int_0^{s}\big(\int_{-sr}^{sr}|\hat{F}(r\hat{\xi},\omega)|^2d\omega\big)r^2\sin\varphi dr d\theta d\varphi\\
&\qquad+\int_0^{2\pi}\int_0^{\pi}\int_{-s}^{0}\big(\int_{-sr}^{sr}|\hat{F}(-r\hat{\xi},\omega)|^2d\omega\big)r^2\sin\varphi dr d\theta d\varphi.
\end{split}
\end{equation}
Let $r=s\hat{r}$ and $\omega=s^2\hat{r}\hat{\omega}$ for $\hat{r},\hat{\omega}\in(-1,1)$. A simple calculation yields 
\begin{equation}
\begin{split}\nonumber
\int_0^{2\pi}&\int_0^{\pi}\int_0^{1}\big(\int_{-s^2\hat{r}}^{s^2\hat{r}}|\hat{F}(s\hat{r}\hat{\xi},\omega)|^2d\omega\big)s^3\hat{r}^2\sin\varphi d\hat{r}d\theta d\varphi\\
&=\int_0^{2\pi}\int_0^{\pi}\int_0^{1}\big(\int_{-1}^{1}|\hat{F}(s\hat{r}\hat{\xi},s^2\hat{r}\hat{\omega})|^2s^2\hat{r}d\hat{\omega}\big)s^3\hat{r}^2\sin\varphi d\hat{r}d\theta d\varphi\\
&=\frac{1}{(2\pi)^4}\int_0^{2\pi}\int_0^{\pi}\int_0^{1}\big( \int_{-1}^1\big|\int_{\R^4}F(x,t)e^{-i(s\hat{r}\hat{\xi}\cdot x+s^2\hat{r}\hat{\omega}t)}dxdt\big|^2s^2\hat{r}d\hat{\omega} \big)s^3\hat{r}^2\sin\varphi d\hat{r}d\theta d\varphi.
\end{split}
\end{equation}
Similarly, we also obtain the expansion of the second term of $(\ref{Ie})$.
This integrals $I(s)$ are analytic functions of $s=s_1+is_2$, $s_1,s_2\in\R$. 
Noting that $|e^{-i(s\hat{r}\hat{\xi}\cdot x+s^2\hat{r}\hat{\omega}t)}|\leq e^{R|s_2|+2T_0|s_1s_2|}$,
we have for all $s\in S$ that 
\begin{equation}\nonumber
|I(s)|\leq CM^2|s|^5e^{2R|s_2|+4T_0|s_1s_2|}.
\end{equation}
Let $\Delta=\max\{2R, 4T_0\}$, it follows for all $s\in S$ that
\begin{equation}\nonumber
|e^{-(\Delta+1)s-(\Delta+1)s^2}I(s)|\leq CM^2.
\end{equation}
Recalling from (\ref{IE}) a prior estimate, we obtain
\begin{equation}\nonumber
|e^{-(\Delta+1)s-(\Delta+1)s^2}I(s)|\leq C\epsilon^2 \quad \mbox{for all} \quad s\in[0,\Lambda].
\end{equation}
A direct application of Lemma \ref{2} shows that 
\begin{equation}\nonumber
|I(s)|\leq CM^2e^{2(\Delta+1)s}\epsilon^{2\mu}, \quad \mbox{for all} s\in(\Lambda, \infty).
\end{equation}
Consider the set (see Figure $\ref{Fig2}$)
\begin{equation}\nonumber
E_1(s)=\{(\xi,\omega)\in \R^4|\ \omega^2-\lambda|\xi|^2=0,\ |\xi|\geq s,\ |\omega|\leq s|\xi|\},
\end{equation}
we have by using the Parseval's identity that 
\begin{equation}\nonumber
\begin{split}
I_1(s):=\int_{E_1(s)}| \hat{F}|^2d\xi\,dw\leq\frac{1}{s^2}\int_{\R^3}|\widehat{\nabla F}|^2d\xi\,dw
\leq\frac{M^2}{s^2}.
\end{split}
\end{equation}
Denote (see Figure $\ref{Fig3}$)
\begin{equation}\nonumber
E_2(s)=\{(\xi,\omega)\in \R^4|\ \omega^2-\lambda|\xi|^2=0,\ \ |\omega|\geq s|\xi|\}.
\end{equation}

Similar to $I(s)$, by using the polar coordinates and $r=s\hat{r}$,  $\hat{r}\in(0,1)$ we get
\begin{equation}\label{i2}
\begin{split}
I_2(s)&=\int_{E_2(s)}|\widehat{F}(\xi,\omega)|^2d\xi d\omega\\
&=\int_0^\infty\big(\int_0^{2\pi}\int_0^{\pi}\int_0^{\frac{1}{s}\omega}|\hat{F}(r\hat{\xi},\omega)|^2r^2\sin\varphi dr\,d\theta\,d\varphi\,\big)d\omega\\
&\qquad+\int_{-\infty}^0\big(\int_0^{2\pi}\int_0^{\pi}\int_0^{-\frac{1}{s}\omega}|\hat{F}(r\hat{\xi},\omega)|^2r^2\sin\varphi dr\,d\theta\,d\varphi\,\big)d\omega.\\
\end{split}
\end{equation}
Simple variable replacement yields
\begin{equation}
\begin{split}
\int_0^\infty&\big(\int_0^{2\pi}\int_0^{\pi}\int_0^{\frac{1}{s}\omega}|\hat{F}(r\hat{\xi},\omega)|^2r^2\sin\varphi dr\,d\theta\,d\varphi\,\big)d\omega\\
&=\int_0^\infty\big(\int_0^{2\pi}\int_0^{\pi}\int_0^{1}|\hat{F}(\frac{1}{s}\omega \hat{r}\hat{\xi},\omega)|^2\frac{1}{s^3}\omega^3 \hat{r}^2\sin\varphi d\hat{r}\,d\theta\,d\varphi\,\big)d\omega\\
&=\frac{1}{s^3}\Big(\int_0^1+\int_1^\infty\big(\int_0^{2\pi}\int_0^{\pi}\int_0^{1}|\hat{F}(\frac{1}{s}\omega \hat{r}\hat{\xi},\omega)|^2\omega^3 \hat{r}^2\sin\varphi d\hat{r}\,d\theta\,d\varphi\,\big)d\omega\Big).
\end{split}
\end{equation}
Similarly, we also obtain the expansion of the second term of $(\ref{i2})$.
Therefore, combining this with
\begin{equation}\nonumber
\begin{split}
\int_0^{1}|\hat{F}(\frac{1}{s}\omega \hat{r}\hat{\xi},\omega)|^2 d\hat{r}
=\int_0^{1} |\frac{1}{(2\pi)^3\omega^3}\int_{\R^3}\partial_t^3F(x,t)e^{-i(\frac{1}{s}\omega \hat{r}\hat{\xi}\cdot x+\omega t)}dxdt|^2d\hat{r},
\end{split}
\end{equation}
we find \[|I_2(s)|\leq C\,\frac{M^2}{s^3}.\]

Now we show the proof of Theorem $\ref{TH3}$.
We assume that $\epsilon<e^{-1}$, otherwise the estimate is obvious. Let 
\begin{eqnarray}
s=
\begin{cases}
\frac{1}{(2(\Delta+2)\pi)^{\frac{1}{4}}}\Lambda^{\frac{1}{2}}|\ln\epsilon|^{\frac{1}{4}(1-\alpha)} & \textrm{if}\ |\ln\epsilon|^{\frac{1}{4}(1-\alpha)}>2^{\frac{1}{4}}\Lambda^{\frac{1}{2}}(2(\Delta+2)\pi)^{\frac{1}{4}},\\
\Lambda & \textrm{if} \ |\ln\epsilon|^{\frac{1}{4}(1-\alpha)}\leq2^{\frac{1}{4}}\Lambda^{\frac{1}{2}}(2(\Delta+2)\pi)^{\frac{1}{3}}.
\end{cases}
\end{eqnarray}
Here $\alpha\in(0,1)$ is given.

Case (i): $|\ln\epsilon|^{\frac{1}{4}(1-\alpha)}>2^{\frac{1}{4}}\Lambda^{\frac{1}{2}}(2(\Delta+2)\pi)^{\frac{1}{4}}$. 
One can check that 
\[s>2^{\frac{1}{4}}\Lambda.\]
Thus, using Lemma \ref{2}, we obtain
\begin{eqnarray}
|I(s)|&\leq & C\, e^{2(\Delta+2)s^2}\epsilon^{2\mu(s)}\cr
&\leq & Ce^{-2\mu(s)|\ln\epsilon|+2(\Delta+2)s^2}\cr
&\leq & Ce^{-\big(\frac{-2|\ln\epsilon|}{\pi}(\frac{\Lambda}{s})^2+2(\Delta+2)s^2\big)}\cr
&\leq & Ce^{-2(\frac{2(\Delta+2)}{\pi})^{\frac{1}{2}}\Lambda|\ln\epsilon|^{1-\frac{1}{2}(1-\alpha)}(1-\frac{1}{2}|\ln\epsilon|^{-\alpha}) }.
\end{eqnarray}
Noting that 
$\frac{1}{2}|\ln\epsilon|^{-\alpha}<\frac{1}{2}$ and $(\frac{2(\Delta+2)}{\pi})^{\frac{1}{2}}>1$,
we have 
\[|I(s)|\leq Ce^{-\Lambda|\ln\epsilon|^{1-\frac{1}{2}(1-\alpha)}}.\]
Using the elementary inequality 
\[e^{-t}\leq \frac{1}{t}, \quad t>0,\]
we get 
\begin{equation}\nonumber
|I(s)|\leq C\frac{M^2}{\big(\Lambda|\ln\epsilon|^{1-\frac{1}{2}(1-\alpha)}\big)}.
\end{equation}
Observing that $\R^4=E(s)\cup E_1(s)\cup E_2(s)$, we know 
\begin{equation}\nonumber
\begin{split}
\|F\|^2_{L^2(\R^4)}=\|\hat{F}\|^2_{L^2(\R^4)}&=\int_{E(s)}|\hat{F}|^2d\xi\,d\omega+\int_{E_1(s)}|\hat{F}|^2d\xi\,d\omega
+\int_{E_2(s)}|\hat{F}|^2d\xi\,d\omega\\
&=I(s)+I_1(s)+I_2(s)\\
&\leq C\big(\frac{M^2}{\Lambda|\ln\epsilon|^{1-\frac{1}{2}(1-\alpha)}}+\frac{M^2}{\Lambda|\ln\epsilon|^{\frac{1}{2}(1-\alpha)}}\big).\\
\end{split}
\end{equation}
Since $\big(\Lambda|\ln\epsilon|^{\frac{1}{2}(1-\alpha)}\big)<\big(\Lambda|\ln\epsilon|^{1-\frac{1}{2}(1-\alpha)}\big)$ when $\Lambda>1$, $|\ln\epsilon|>1$ and $0<\alpha<1$, we obtain 
\begin{equation}\label{F1}
\|F\|^2_{L^2(\R^4)}\leq C\frac{M^2}{\Lambda|\ln\epsilon|^{\frac{1}{2}(1-\alpha)}}. 
\end{equation}

Case (ii): $|\ln\epsilon|^{\frac{1}{4}(1-\alpha)}\leq2^{\frac{1}{4}}\Lambda^{\frac{1}{2}}(2(\Delta+2)\pi)^{\frac{1}{4}}$. In this case we have for $s=\Lambda$ that
\[|I(s)|=|I(\Lambda)| \leq |E(\Lambda)|\Lambda^4\,(1+\Lambda)^2\epsilon^2.\]
Using  estimates of $I_1(s)$ and $I_2(s)$ , we obtain 
\begin{equation}\label{Hest0}
\begin{split}
\|F\|^2_{L^2(\R^4)}=&\|\hat{F}\|^2_{L^2(\R^4)}\\
=&I(s)+I_1(s)+I_2(s)\\
\leq& C\Big(|E(\Lambda)|\Lambda^4\,(1+\Lambda)^2\epsilon^2+\frac{M^2}{\Lambda|\ln\epsilon|^{\frac{1}{2}(1-\alpha)}}\Big). 
\end{split}
\end{equation}
Combining $(\ref{F1})$ and $(\ref{Hest0})$, we finally derive 
\begin{equation}\nonumber
\|F\|^2_{L^2(\R^4)}
\leq C\Big(\Lambda^{10}\,\epsilon^2 +\frac{M^2}{\Lambda|\ln\epsilon|^{\frac{1}{2}(1-\alpha)}}\Big). 
\end{equation}
This completes the proof.

\section{Proof of Theorem \ref{TH4}}\label{th3}
In this section, we assume that $F$ takes the form
\begin{equation}\label{3F0}
F(x_1,x_2,x_3,t)=f(\tilde{x},t)g(x_3), \quad \tilde{x}=(x_1, x_2)\in\R^2,  \quad   t\in(0,+\infty).
\end{equation}
%where $f$ is compactly supported on $\widetilde{B}_{R_0}\times(0,T_0)$ for some $0<R_0<R/\sqrt{2}$ and $g\in H^1(\R)$ is supported on $(-R_0,R_0)$. Here $\tilde{x}=(x_1,x_2)\in\R^2$ for $x=(x_1,x_2,x_3)\in\R^3$ and $\widetilde{B}_{R_0}:=\{\tilde{x}\in\R^2:\ |\tilde{x}|<R_0\}$.
Then the equation $(\ref{eq1})$ becomes
\begin{equation}\label{eqn5}
\partial_t^2u-\Delta u=f(\tilde{x},t)\,g(x_3).
\end{equation}
Assuming that $g$ is known, we establish an increasing stability estimate for $f$ from the Dirichlet data $\{u(x,t)| \ x\in\partial B_R, \ t\in(0,T)\}$.
We choose the test function
$$w(\xi, t)=e^{-i(\xi_1\cdot x_1+\xi_2\cdot x_2+\xi_3x_3+\omega t)},\quad
\omega^2-|\xi_1|^2-|\xi_2|^2=|\xi_3|^2,
$$ which satisfies the wave equation
\begin{equation}\nonumber
\partial_t^2w-\Delta w=0.
\end{equation}
By the  Huygens' principle, it holds that $u(x,t)=0$ for $|x|<R$ and $t>T$.
Multiplying $w$ on both sides of $(\ref{eqn5})$ and integrating over $B_R\times(0,T)$, we have
\begin{equation}\label{fg}
\begin{split}
(2\pi)^{2}\hat{f}(\xi_1,\xi_2,\omega)\hat{g}(\xi_3)
&=\int_{\R^3}f(\tilde{x},t)e^{-i(\tilde{\xi}\cdot \tilde{x}+\omega t)}d\tilde{x}dt\int_{\R}g(x_3)e^{-i \xi_3\cdot x_3}dx_3\\
&=\int_0^T\int_{\partial B_R}(u\frac{\partial w}{\partial\nu}-w\frac{\partial u}{\partial\nu})ds(x)dt,
\end{split}
\end{equation}
where $\tilde{\xi}=(\xi_1,\xi_2)$.
Consider the set $$E(b)=\{(\xi_1,\xi_2,\omega)\in\R^3| \ \omega^2-|\xi_1|^2-|\xi_2|^2=|\xi_3|^2, \ \xi_3\in(-b,b), \ |\xi_1|^2+|\xi_2|^2\leq b^2\}.$$
Since $\xi_3\in(-b,b)$, one can check that $|E_b|>0$.
Using the fact that $f$ and $g$ have compact supports and arguing analogously to Lemma \ref{lem4}, we deduce
\begin{equation}\label{6u}
 \|\partial_{\nu}u\|_{L^2([0,T];H^{\frac{1}{2}}(\partial B_R))}
\leqslant C\|u\|_{H^2([0,T];H^{\frac{3}{2}}(\partial B_R))}.
\end{equation}
Combining (\ref{6u}) and  $|\hat{g}(\xi_3)|\geq\delta>0$ for $\xi_3\in(-b,b)$ with $(\ref{fg})$, we get for $(\xi_1,\xi_2,\omega)\in E(b)$ that
\begin{equation}\label{3hf}
 |\hat{f}(\xi_1,\xi_2,\omega)|\leq C\frac{(1+|\xi|)\|u\|_{H^2([0,T];H^\frac{3}{2}(\partial B_R))}}{|\hat{g}(\xi_3)|}\leq C \delta^{-1} (1+|\xi|) \epsilon,
\end{equation}
where $|\xi|^2=\xi_1^2+\xi_2^2+\xi_3^2$.

\begin{figure}[htbp]
 \centering	
\includegraphics[width=0.48\textwidth]{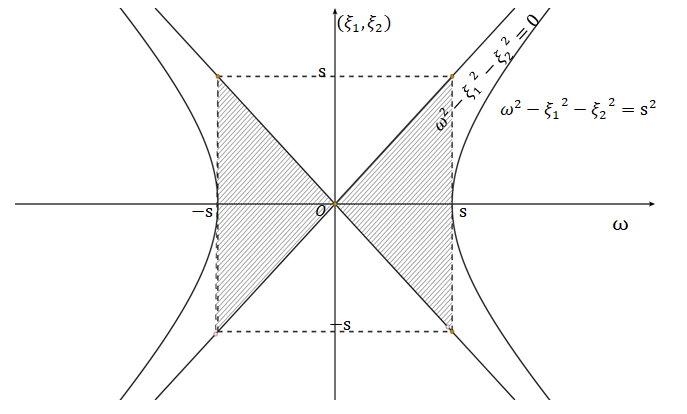}
	\caption{$E(s)$ is the shaded area.}
	\label{Fig4}
\end{figure}

\begin{figure}[htbp]
 \centering	
\includegraphics[width=0.48\textwidth]{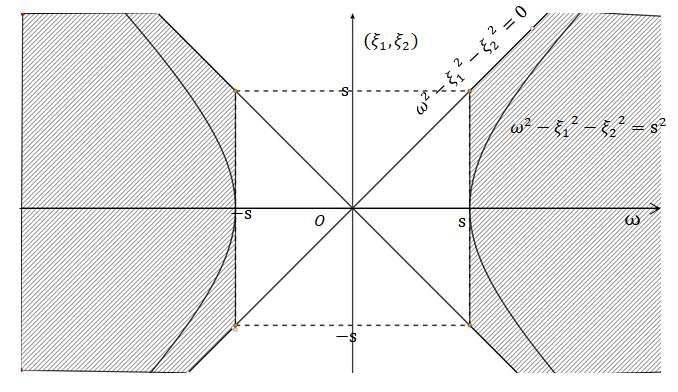}
	\caption{$E_1(s)$ is the shaded area.}
	\label{Fig5}
\end{figure}

\begin{figure}[htbp]
 \centering	
\includegraphics[width=0.48\textwidth]{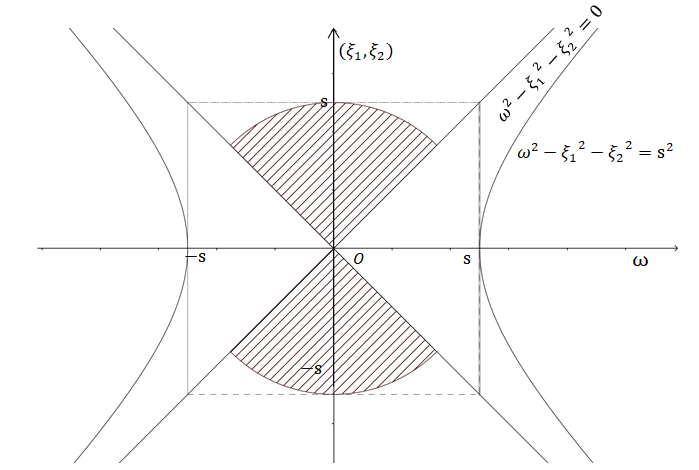}
	\caption{$E_2(s)$ is the shaded area.}
	\label{Fig6}
\end{figure}

Define the set (see Figure $\ref{Fig4}$)
\begin{equation}\nonumber
E(s)=\{(\xi_1,\xi_2,\omega)\in\R^3| \ \omega^2-|\xi_1|^2-|\xi_2|^2=|\xi_3|^2, \ \xi_3\in(-s,s), \ |\xi_1|^2+|\xi_2|^2\leq s^2\}.
\end{equation}
Let
$$I(s)=\int_{E(s)}|\hat{f}(\xi_1,\xi_2,\omega)|^2d\xi_1d\xi_2d\omega.$$ 
Using the polar coordinates $\xi_1=r\sin\theta$,  $\xi_2=r\cos\theta$,
$0\leq r\leq s$,  $0\leq\theta\leq2\pi$, we obtain that
\begin{equation}\nonumber
\begin{split}
I(s)=\int_{-s}^0\int_0^{2\pi}\int_0^{-\omega}|\hat{f}(r\sin\theta,r\cos\theta,\omega)|^2 rdrd\theta d\omega+\int_{0}^s\int_0^{2\pi}\int_0^{\omega}|\hat{f}(r\sin\theta,r\cos\theta,\omega)|^2 rdrd\theta d\omega.
\end{split}
\end{equation}
Let $\omega=s\overline{\omega}$ for $\overline{\omega}\in (-1,1)$ and $r=s\overline{r}$ for $\overline{r}\in (0,1)$. Then a simple calculation yields 
\begin{equation}\nonumber
\begin{split}
I(s)=\int_{-1}^0\int_0^{2\pi}&\int_0^{-\overline{\omega}}|\hat{f}(s\overline{r}\sin\theta,s\overline{r}\cos\theta,s\overline{\omega})|^2s^3\overline{r}d\overline{r}d\theta d \overline{\omega}\\
&\qquad+\int_{0}^1\int_0^{2\pi}\int_0^{\overline{\omega}}|\hat{f}(s\overline{r}\sin\theta,s\overline{r}\cos\theta,s\overline{\omega})|^2s^3\overline{r}d\overline{r}d\theta d \overline{\omega}.\\
\end{split}
\end{equation}
This integral $I(s)$ is also an analytic function of $s=s_1+is_2\in\C$.
Noting that 
$|e^{-i(s\overline{r}(sin \theta,cos \theta)\cdot\tilde{x}+s\overline{\omega}t)}|\leq e^{(R_0+T_0)|s_2|}$ for all $\tilde{x}\in \tilde{B}_{R_0}$ and $t\in (0,T_0)$, we deduce 
\begin{equation}
\begin{split}
(2\pi)^{\frac{3}{2}}|\hat{f}(s\overline{r}\sin\theta,s\overline{r}\cos\theta,s\overline{\omega})|^2&=|\int_0^{T_0}\int_{\tilde{B}_{R_0}} {f}(\tilde{x},t)e^{-i(s\overline{r}(sin \theta,cos \theta)\cdot\tilde{x}+s\overline{\omega}t)} d \overline{x}d t|^2\\
&\leq Ce^{2(R_0+T_0)|s_2|}\|f\|_{L^2(\mathbb{R}^3)}^2.
\end{split}
\end{equation}
Thus
\begin{equation}\nonumber
|I(s)|\leq C|s|^3e^{2(R_0+T_0)|s_2|}\|f\|_{L^2(\mathbb{R}^3)}^2,
\end{equation} 
which yields
\begin{equation}\nonumber
|I(s)e^{-\big(2(R_0+T_0)+1\big)s}|\leq CM^2\  \mbox{for all} \ s\in S.
\end{equation} 
Together with the estimate (\ref{3hf}) and definition of $E(s)$, we have 
\begin{equation}\label{3Cb}
\begin{split}
|I(s)|\leq C\,|E(s)|\,(1+s)^2\,\epsilon^2 \quad \mbox{for all} \quad s\in[0,b].
\end{split}
\end{equation}
Applying Lemma \ref{2}, we know for all $s>b$ that 
\begin{equation}\nonumber
\begin{split}
|I(s)|\leq CM^2e^{\big(2(R_0+T_0)+1\big)s}\epsilon^{2\mu(s)}.
\end{split}
\end{equation}
Define $E_1(s):=\{(\xi_1,\xi_2,\omega)\in\R^3| \, |\omega|>s, \ \xi_1^2+\xi_2^2 \leq\omega^2 \}$ (see Figure $\ref{Fig5}$) and 
\[I_1(s)=\int_{E_1(s)}|\hat{f}(\xi_1,\xi_2,\omega)|^2d\xi_1d\xi_2d\omega.\]
Since supp$f(\tilde{x},t)\subset \widetilde{B}_{R_0}\times(0,T_0)$, we have
\begin{equation}
\nonumber
\begin{split}
(2\pi)^{\frac{3}{2}}\hat{f}(\xi_1,\xi_2,\omega)&=\int_0^{T_0}\int_{\tilde{B}_{R_0}} {f}(x_1,x_2,t)e^{-i(\xi_1\cdot x_1+\xi_2\cdot x_2+\omega t)}d x_1d x_2d t\\
&=\dfrac{1}{-i\omega}\int_0^{T_0}\int_{\tilde{B}_{R_0}} \partial_t f (x_1,x_2,t)e^{-i(\xi_1\cdot x_1+\xi_2\cdot x_2+\omega t)}d x_1d x_2d t\\
&=\dfrac{1}{-i\omega}\widehat{\partial_t f}.
\end{split}
\end{equation}
Then using the Parseval's identity, it's easy to get
\begin{equation}
I_1(s)=\int_{E_1(s)}|\widehat{f}(\xi_1,\xi_2,\omega)|^2d\xi_1d\xi_2d\omega\leq \frac{1}{s^2}\int_{E_1(s)}|\widehat{ \nabla f}(\xi_1,\xi_2,\omega)|^2d\xi_1d\xi_2d\omega \leq\frac{CM^2}{s^2}.
\end{equation}
Define the set (see Figure $\ref{Fig6}$)
\begin{equation}\nonumber
E_2(s):=\{(\xi_1,\xi_2,\omega)\in\R^3|\, \omega^2+|\xi_1|^2+|\xi_2|^2\leq s^2, \ \xi_1^2+\xi_2^2\geq\omega^2 \}.
\end{equation}

\begin{figure}[htbp]
 \centering	
\includegraphics[width=0.48\textwidth]{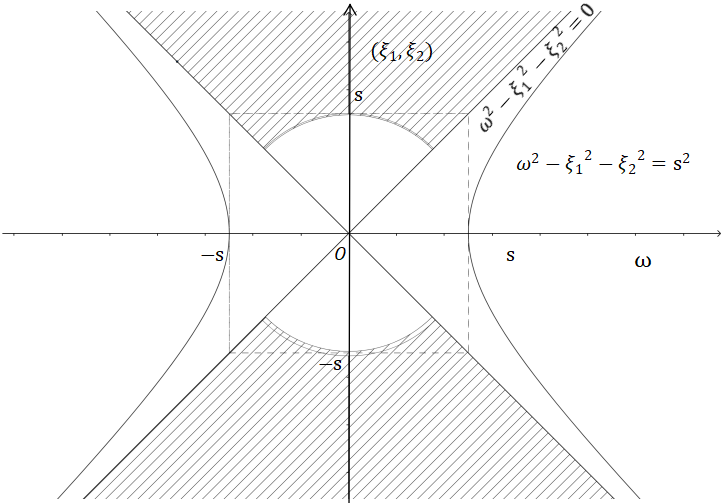}
	\caption{$E_3(s)$ is the shaded area.}
	\label{Fig7}
\end{figure}

Now we estimate  
$$I_2(s)=\int_{E_2(s)}|\hat{f}(\xi_1,\xi_2,\omega)|^2d\xi_1d\xi_2d\omega.$$
Define
$F_b(\xi_1,\xi_2,\omega)=\hat{f}(b\xi_1,b\xi_2,b\omega)$ for any $(\xi_1,\xi_2,\omega)\in\mathbb{R}^3.$  Since $f$ is compactly supported, one can see that the function $F_b$ is analytic and it satisfies for $\gamma\in (\N\cup\{0\})^3$ that
\begin{eqnarray*}
|\partial^\gamma F_b(\xi_1,\xi_2,\omega)|&=&|\partial^\gamma F(b\xi_1,b\xi_2,b\omega)|\\
&=&\Big|\partial^\gamma \int_{\mathbb{R}^3}f(x_1,x_2,t) e^{-ib(\xi_1\cdot x_1+\xi_2\cdot x_2+\omega t)}\,dx_1 dx_2dt\Big|\cr
&=& \Big|  \sum_{\alpha+\beta=\gamma \atop \alpha,\beta\in\N\cup0}\int_{\mathbb{R}^3} (-i)^{|\gamma|} b^{|\gamma|} x^{\alpha}t^{\beta} f(x_1,x_2,t) e^{-ib(\xi_1\cdot x_1+\xi_2\cdot x_2+\omega t)}\,dx_1 dx_2dt \Big|.
\end{eqnarray*}
Using $b^{|\gamma|}<|\gamma|!\, e^b $ and  $\eta=\max\{R,T\}^{-1}$, 
we obtain
\begin{eqnarray}
|\partial^\gamma F_b(\xi_1,\xi_2,\omega)|\leq \|f\|^2_{L^2(B_R)}\, \eta^{-|\gamma|} \,b^{|\gamma|} 
\leq  C\|f\|^2_{L^2(B_R)}\,\eta^{-|\gamma|}\, |\gamma|!\, e^b \leq  C M^2\,\eta^{-|\gamma|}\, |\gamma|!\, e^b.
\end{eqnarray}
Applying Proposition \ref{pro3.1} to the set $\mathcal{O}$ defined as $\mathcal{O}:= E(1)$,
we can find a constant $\alpha\in (0,1)$ such that
 $$\|F_b\|_{L^\infty(B(0,1))}\leq C\, e^{b(1-\alpha)} \|F_b\|^\alpha_{L^\infty(\mathcal{O})}.$$
Using the fact that $\hat{f}(\xi)=F_b(b^{-1}\xi)$, one gets the following estimate
\begin{eqnarray}\nonumber
\|\hat{f}\|_{L^{\infty}(B(0,b))}&=&\|F_b\|_{L^\infty(B(0,1))}\leq Ce^{b(1-\alpha)} \|F_b\|_{L^\infty(\mathcal{O})}^\alpha\\
&\leq& C e^{b(1-\alpha)} \|\hat{f}\|_{L^\infty(E(b))}^\alpha
\leq Ce^{b(1-\alpha)} (1+b)^\alpha\epsilon^\alpha.
\end{eqnarray}
Combining (\ref{3hf}), we know for $s\in[0,b]$ that 
\begin{equation}\nonumber
|I_2(s)|\leq C|E_2(b)|\,e^{2b(1-\alpha)}\, (1+b)^{2\alpha}\,\epsilon^{2\alpha}
\end{equation}
and similarly 
\begin{equation}\nonumber
\begin{split}
|I_2(s)|\leq CM^2\,e^{\big(2(R_0+T_0)+1\big)|s_2|} \quad \mbox{for all} \quad s\in S.
\end{split}
\end{equation}
Thus applying Lemma \ref{2}, we have 
\begin{equation}\nonumber
\begin{split}
|e^{-[2(R_0+T_0)+1]s}I_2(s)|\leq CM^2(\epsilon^{\alpha})^{2\mu} \quad \mbox{for all} \quad s>b.
\end{split}
\end{equation}
Consider the set (see Figure $\ref{Fig7}$)
\begin{equation}\nonumber
E_3(s):=\{(\xi_1,\xi_2,\omega)\in\R^3| \, \ \omega^2+|\xi_1|^2+|\xi_2|^2> s^2, \ \xi_1^2+\xi_2^2\geq\omega^2 \}.
\end{equation}

As done in the estimate of $I_1(s)$, appying Parseval's identity yields
\begin{equation}
I_3(s)=\int_{E_3(s)}|\hat{f}(\xi_1,\xi_2,\omega)|^2d\xi_1d\xi_2d\omega\leq\frac{CM^2}{s^2}.
\end{equation}
Now we prove Theorem \ref{TH4}.
We assume that $\epsilon<e^{-1}$, since if otherwise the estimate is obvious. Let 
\begin{eqnarray}
s=
\begin{cases}
\frac{1}{((2(R_0+T_0)+3)\pi)^{\frac{1}{3}}}b^{\frac{2}{3}}|\ln\epsilon|^{\frac{1}{4}} & \textrm{if}\ |\ln\epsilon|^{\frac{1}{4}}>2^{\frac{1}{4}}b^{\frac{1}{3}}((2(R_0+T_0)+3)\pi)^{\frac{1}{3}},\\
b & \textrm{if} \ |\ln\epsilon|^{\frac{1}{4}}\leq 2^{\frac{1}{4}}b^{\frac{1}{3}}((2(R_0+T_0)+3)\pi)^{\frac{1}{3}}.
\end{cases}
\end{eqnarray}

Case (i): $|\ln\epsilon|^{\frac{1}{4}}>2^{\frac{1}{4}}b^{\frac{1}{3}}((2(R_0+T_0)+3)\pi)^{\frac{1}{3}}$. 
One can check that 
\[s>2^{\frac{1}{4}}b.\]
Thus, using Lemma \ref{2} we obtain
\begin{eqnarray}
|I(s)|&\leq & CM^2\, e^{(2(R_0+T_0)+3)s}\epsilon^{2\mu}\cr
&\leq & Ce^{-2\mu|\ln\epsilon|+((2(R_0+T_0)+3)s}\cr
&\leq & CM^2\,e^{-\frac{-2|\ln\epsilon|}{\pi}(\frac{b}{s})^2+\frac{(2(R_0+T_0)+3)}{((2(R_0+T_0)+3)\pi)^{\frac{1}{3}}}b^{\frac{2}{3}}|\ln\epsilon|^{\frac{1}{4}}}\cr
&\leq & CM^2\, e^{-2\big(\frac{(2(R_0+T_0)+3)^2}{\pi}\big)^{\frac{1}{3}}b^{\frac{2}{3}}|\ln\epsilon|^{\frac{1}{2}}\big(1-\frac{1}{2} |\ln\epsilon|^{-\frac{1}{4}} \big)}.
\end{eqnarray}
Noting that 
$\frac{1}{2}|\ln\epsilon|^{-\frac{1}{4}}<\frac{1}{2}$ and $\big(\frac{(2(R_0+T_0)+3)^2}{\pi}\big)^{\frac{1}{3}}>1$,
we have 
\[|I(s)|\leq C M^2\,e^{-b^{\frac{2}{3}}|\ln\epsilon|^{\frac{1}{2}}}.\]
Using the elementary inequality 
\[e^{-t}\leq \frac{C}{t^2}, \quad t>0,\]
we get 
\begin{equation}\nonumber
|I(s)|\leq\frac{CM^2}{\big(b^{\frac{2}{3}}|\ln\epsilon|^{\frac{1}{2}}\big)^2}.
\end{equation}
Similarly, we have
\begin{equation}\nonumber
|I_2(s)|\leq \frac{CM^2}{\big(b^{\frac{2}{3}}|\alpha\ln\epsilon|^{\frac{1}{2}}\big)^2}.
\end{equation}
Thus
\begin{equation}\nonumber
\begin{split}
\|f\|^2_{L^2(\R^3)}&=\|\hat{f}\|^2_{L^2(\R^3)}\\
&=\int_{E(s)}|\hat{f}(\tilde{\xi},\omega)|^2\,d\tilde{\xi}\,d\omega+\int_{E_1(s)}|\hat{f}(\tilde{\xi},\omega)|^2\,d\tilde{\xi}\,d\omega\\
&\qquad+\int_{E_2(s)}|\hat{f}(\tilde{\xi},\omega)|^2\,d\tilde{\xi}\,d\omega+\int_{E_3(s)}|\hat{f}(\tilde{\xi},\omega)|^2\,d\tilde{\xi}\,d\omega\\
&=I(s)+I_1(s)+I_2(s)+I_3(s)\\
&\leq C\Big(\frac{M^2}{\big(b^{\frac{2}{3}}|\ln\epsilon|^{\frac{1}{2}}\big)^2}+
\frac{M^2}{\big(b^{\frac{2}{3}}|\ln\epsilon|^{\frac{1}{4}}\big)^2}+\frac{M^2}{\big(b^{\frac{2}{3}}|\alpha\ln\epsilon|^{\frac{1}{2}}\big)^2}+\frac{M^2}{\big(b^{\frac{2}{3}}|\alpha\ln\epsilon|^{\frac{1}{4}}\big)^2}
\Big).
\end{split}.
\end{equation}
Since $b^{\frac{2}{3}}|\ln\epsilon|^{\frac{1}{2}}>b^{\frac{2}{3}}|\ln\epsilon|^{\frac{1}{4}}$, $b^{\frac{2}{3}}|\alpha\ln\epsilon|^{\frac{1}{2}}>b^{\frac{2}{3}}|\alpha\ln\epsilon|^{\frac{1}{4}}$ and $b^{\frac{2}{3}}|\ln\epsilon|^{\frac{1}{4}}>b^{\frac{2}{3}}|\alpha\ln\epsilon|^{\frac{1}{4}}$ when $b>1$ and $|\ln\epsilon|>1$, we obtain
\begin{equation}\label{I1}
\|f\|^2_{L^2(\R^3)}\leq C\frac{M^2}{\big(b^{\frac{2}{3}}|\alpha\ln\epsilon|^{\frac{1}{4}}\big)^2}.
\end{equation}

Case (ii): $|\ln\epsilon|^{\frac{1}{4}}>2^{\frac{1}{4}}b^{\frac{1}{3}}((2(R_0+T_0)+3)\pi)^{\frac{1}{3}}$. In this case we have from  $s=b$ and (\ref{3Cb}) that
\[|I(s)|=|I(b)| \leq |E(b)|(1+b)^2\epsilon^2.\]
Combining this estimate and $I_1(s)$, we obtain 
\begin{equation}\label{I2}
\begin{split}
\int_{E(s)}|\hat{f}(\tilde{\xi},\omega)|^2\,d\tilde{\xi}\,d\omega &+\int_{E_1(s)}|\hat{f}(\tilde{\xi},\omega)|^2\,d\tilde{\xi}\,d\omega
=I(s)+I_1(s)\\
&\leq C\Big(b^5\epsilon^2
+\frac{M^2}{\big(b^{\frac{2}{3}}|\ln\epsilon|^{\frac{1}{4}}\big)^2}\Big).
\end{split}
\end{equation}
Similarly, we have 
\begin{equation}\nonumber
|I_2(s)|=|I_2(b)|\leq |E_2(b)|e^{2b(1-\alpha)} (1+b)^{2\alpha}\epsilon^{2\alpha}
\end{equation}
and
\begin{equation}\label{I3}
\begin{split}
\int_{E_2(s)}|\hat{f}(\tilde{\xi},\omega)|^2\,d\tilde{\xi}\,d\omega &+\int_{E_3(s)}|\hat{f}(\tilde{\xi},\omega)|^2\,d\tilde{\xi}\,d\omega
=I_2(s)+I_3(s)\\
&\leq C\Big(b^3\,\,e^{2b(1-\alpha)}(1+b)^{2\alpha}\epsilon^{2\alpha}+\frac{M^2}{\big(b^{\frac{2}{3}}|\alpha\ln\epsilon|^{\frac{1}{4}}\big)^2}\Big).
\end{split}
\end{equation}
Combining $(\ref{I1})$, $(\ref{I2})$ and $(\ref{I3})$, we obtain 
\begin{equation}\nonumber
\begin{split}
\|f\|^2_{L^2(\R^3)}&=\|\widehat{f}\|^2_{L^2(\R^3)}\\
&=I(s)+I_1(s)+I_2(s)+I_3(s)\\
&\leq C\Big(b^5\epsilon^2+b^3\,\,e^{2b(1-\alpha)}(1+b)^{2\alpha}\epsilon^{2\alpha}
+\frac{M^2}{\big(b^{\frac{2}{3}}|\alpha\ln\epsilon|^{\frac{1}{4}}\big)^2}\Big).
\end{split}
\end{equation}
This completes the proof.

\section{Conclusion}\label{th4}
In this work, we have shown increasing stability estimates of the $L^2$-norm of the source function with measurement data taken on a sphere over a long time interval in $\R^3$. These results are expected to be valid in the two-dimensional case and also for  more general evolutional equations. The absence of Huygens' principle in 2D will create  technical difficulties and a new technique must be developed. We hope to be able to report the progress on these problems in the future.

\section{Acknowledgment}
The work of G. Hu is partially supported by the National Natural Science Foundation of China (No. 12071236) and the Fundamental Research Funds for Central Universities in China (No. 63233071). The work of S. Si is supported by the Natural Science Foundation of Shandong Province, China (No. ZR202111240173).

%\appendix
% \section{Two theorems}
% We list two important theorems used in this paper. The following is Tichmarsh convolution theorem \cite{titchmarsh}.
% \begin{theorem}\label{A1}
% If $\varphi(t)$ and $\psi(t)$ are integrable functions, such that
% \begin{equation}
% \nonumber
% \int_0^x \varphi(t)\psi(x-t)\mathrm{d} t=0
% \end{equation}
% almost everywhere in interval $0<x<\kappa$, then $\varphi(t)=0$ almost everywhere in $(0,\lambda)$, and $\psi(t)=0$
% almost everywhere in $(0,\mu)$, where $\lambda+\mu\geqslant\kappa$.
% \end{theorem}
%
% Next is Firtz John's global Holmgren theorem  (\cite[\S 3.5]{Bers1964}, \cite[\S 1.7]{Rauch1991}).
%Assume that $u$ is a local solution of the $\mathit{m}$th order linear Cauchy problem
%\begin{align}\label{eqA}
%\begin{cases}
%  \mathit{P}u=f   &\textrm{in} \ \Omega,\\
%\partial^{\alpha}u=0  &\textrm{on} \ \Sigma \ \textrm{for all}\  |\alpha|\leq \mathit{m}-1.
%\end{cases}
%\end{align}
%\begin{theorem}
%Suppose that $\mathit{P}$ is a linear partial differential operator with coefficients real analytic on a neighborhood of
%$\tilde{x}\in \mathbb{R}^3 $, and $\Sigma$ is a $C^{\mathit{m}}$ embedded hypersurface noncharacteristic at $\tilde{x}$. If $u$ is a $C^{\mathit{m}}$ solution of $(\ref{eqA})$ on a neighborhood of $\tilde{x}$, then $u$ vanishes on a neighborhood of $\tilde{x}$.
%\end{theorem}
%In this paper, we take the surface $\Sigma$ as the sphere, which is obviously noncharacteristic. And the coefficients are all constants, which means that the theorem works.

\end{document}